\documentclass[a4paper]{amsart}

\usepackage[dvipsnames]{xcolor} 
\usepackage{amsmath,amsthm,amssymb}
\usepackage{enumitem}
\usepackage{tikz}
\usepackage{hyperref}
\usepackage{amsfonts,marvosym,amscd}
\usepackage{mathtools}
\usepackage{mathscinet}
\usepackage{wasysym} 
\usepackage{graphicx}
\usepackage{psfrag}
\usepackage{datetime}
\usepackage{tikz-cd} 
\usepackage{url}
\usepackage{adjustbox}
\usepackage{pifont} 
\usepackage{ytableau}
\usepackage{bbm} 

\usepackage{tabstackengine}
\setstackgap{L}{.75\baselineskip}

\usepackage{pdflscape}
\usepackage{afterpage}
\usepackage{capt-of}

\usetikzlibrary{decorations.pathmorphing}
\usetikzlibrary{decorations.pathreplacing}
\usetikzlibrary{intersections}
\usetikzlibrary{shapes,backgrounds,shapes.misc, positioning,shapes.geometric,arrows,matrix,fit,calc}

\makeatletter
\tikzset{
    dot diameter/.store in=\dot@diameter,
    dot diameter=3pt,
    dot spacing/.store in=\dot@spacing,
    dot spacing=10pt,
    dots/.style={
        line width=\dot@diameter,
        line cap=round,
        dash pattern=on 0pt off \dot@spacing
    }
}
\makeatother

\pgfdeclarelayer{foreground} 
\pgfdeclarelayer{background}
   \pgfsetlayers{background,%
                 main,%
                 foreground%
                 }

\numberwithin{equation}{section}

\newtheorem{theorem}{Theorem}[section]
\newtheorem{lemma}[theorem]{Lemma}
\newtheorem{proposition}[theorem]{Proposition}

\newtheorem{conjecture}[theorem]{Conjecture}

\theoremstyle{definition}
\newtheorem{definition}[theorem]{Definition}

\newtheorem{example}[theorem]{Example}

\newtheorem{construction}[theorem]{Construction}

\theoremstyle{remark}
\newtheorem{remark}[theorem]{Remark}

\newcommand{\cmark}{\ding{51}}
\newcommand{\xmark}{\ding{55}}

\let\emptyset\varnothing

\newcommand{\spt}[1]{\{\,#1\,\}} 
\newcommand{\st}{\mid}

\DeclareFontFamily{U}{mathx}{\hyphenchar\font45}
\DeclareFontShape{U}{mathx}{m}{n}{
      <5> <6> <7> <8> <9> <10>
      <10.95> <12> <14.4> <17.28> <20.74> <24.88>
      mathx10
      }{}
\DeclareSymbolFont{mathx}{U}{mathx}{m}{n}
\DeclareFontSubstitution{U}{mathx}{m}{n}
\DeclareMathAccent{\widecheck}{0}{mathx}{"71}

\newcommand{\equ}{\theta} 
\newcommand{\qrel}{\mathrel{\leqslant_{\equ}}} 
\newcommand{\set}{S} 
\newcommand{\rel}{R} 
\newcommand{\quot}{P/{\equ}} 
\newcommand{\lquot}{L/{\equ}} 
\newcommand{\dual}[1]{\overline{#1}} 

\newcommand{\meet}{\wedge}
\newcommand{\join}{\vee}

\newcommand{\bruhat}[2]{\mathcal{B}(#1, #2)} 
\newcommand{\stash}[2]{\mathcal{S}(#1, #2)} 

\newcommand{\std}[1]{\mathrm{std}(#1)} 

\newcommand{\collapse}[1]{\mathsf{coll}(#1)} 
\newcommand{\pleq}{\preccurlyeq} 
\newcommand{\pgeq}{\succcurlyeq} 

\newcommand{\posetcat}{\mathbf{Poset}} 
\newcommand{\preordcat}{\mathbf{Pre}\text{-}\mathbf{ord}} 

\DeclareMathOperator{\Hom}{Hom}

\newcommand{\ls}[2]{L_{#1}(#2)} 
\newcommand{\us}[2]{U_{#1}(#2)} 
\newcommand{\mac}[1]{\mathsf{L}(#1)} 
\newcommand{\maco}[1]{\mathsf{L}_{0}(#1)} 
\newcommand{\macom}[1]{\iota_{#1}} 
\newcommand{\lex}[1]{\mathsf{lex}(#1)}

\newcommand{\indic}[2]{\mathbbm{1}_{#1}(#2)} 

\newcommand{\trans}[1]{\overrightarrow{#1}}

\DeclareMathOperator{\conv}{conv} 
\newcommand{\ip}[2]{\langle #1, #2 \rangle} 
\DeclareMathOperator{\argmax}{argmax}

\usepackage[backend=biber,
    giveninits=true,
	style=alphabetic,
	url=false,
	doi=false,
	isbn=false,
	maxbibnames=99]{biblatex}
\DeclareFieldFormat{title}{#1}
\DeclareFieldFormat{postnote}{#1} 
\DeclareFieldFormat{volume}{Volume #1}
\DeclareFieldFormat{edition}{#1 edition}
\renewbibmacro{in:}{} 

\makeatletter
\def\@tocline#1#2#3#4#5#6#7{\relax
  \ifnum #1>\c@tocdepth 
  \else
    \par \addpenalty\@secpenalty\addvspace{#2}%
    \begingroup \hyphenpenalty\@M
    \@ifempty{#4}{%
      \@tempdima\csname r@tocindent\number#1\endcsname\relax
    }{%
      \@tempdima#4\relax
    }%
    \parindent\z@ \leftskip#3\relax \advance\leftskip\@tempdima\relax
    \rightskip\@pnumwidth plus4em \parfillskip-\@pnumwidth
    #5\leavevmode\hskip-\@tempdima
      \ifcase #1
       \or\or \hskip 1em \or \hskip 2em \else \hskip 3em \fi%
      #6\nobreak\relax
    \dotfill\hbox to\@pnumwidth{\@tocpagenum{#7}}\par
    \nobreak
    \endgroup
  \fi}
\makeatother

\title[Congruences and quotients of posets]{A survey of congruences and quotients of partially ordered sets}
\author{Nicholas J. Williams}
\urladdr{https://nchlswllms.github.io/}
\email{nw480@cam.ac.uk}
\address{Pavilion E, Department of Pure Mathematics and Mathematical Statistics, Centre for Mathematical Sciences, University of Cambridge, Wilberforce Road, Cambridge, CB3 0WB, United Kingdom}
\thanks{This work was partially supported by~EPSRC grant EP/V050524/1 and by a JSPS International Short-Term Postdoctoral Research Fellowship at the University of Tokyo. I would like to thank Mikhail Gorsky and an anonymous referee for helpful feedback on this survey and Russ Woodroofe for providing links to some useful references.}

\begin{document}

\begin{abstract}
A quotient of a poset $P$ is a partial order obtained on the equivalence classes of an equivalence relation $\equ$ on $P$; $\equ$ is then called a congruence if it satisfies certain conditions, which vary according to different theories. The literature on congruences and quotients of partially ordered sets contains a large and profilerating array of approaches, but little in the way of systematic exposition and examination of the subject. We seek to rectify this by surveying the different theories in the literature and providing philosophical discussion on requirements for notions of congruences of posets. We advocate a pluralist approach which recognises that different types of congruence arise naturally in different mathematical situations. There are some notions of congruence which are very general, whilst others capture specific structure which often appears in examples. Indeed, we finish by giving several examples where quotients of posets appear naturally in mathematics.
\end{abstract}

\maketitle

\tableofcontents

\section{Introduction}

Quotients and partially ordered sets are among the most basic notions in mathematics, and yet their interaction with each other has received little systematic study. Roughly speaking, a quotient of a poset $P$ is a poset $Q$ whose elements are equivalence classes of an equivalence relation on $P$, and whose order relations are determined by those of $P$ in a natural way. Another way of viewing $Q$ is as a certain poset with a surjective order-preserving map $P \to Q$; here the fibres of the map correspond to the equivalence relation. Given an arbitrary equivalence relation on~$P$, there is not always a natural way of constructing a poset whose elements are the equivalence classes. Hence, one needs to restrict the class of equivalence relations, or take a different approach to quotients.
Such a restricted class of equivalence relations is known as a class of \emph{congruences}.

There is clear motivation for the notion of a quotient of a poset, and therefore of a congruence. Given a set which is partially ordered and partitioned, it is natural to ask whether the partial order can be used to give a natural partial order to the parts of the partition. Furthermore, quotients allow one to construct new posets from old ones, and are useful for describing the relation between existing posets.

Several different approaches to congruences of posets exist \cite{sturm_ao,stanley_quotient_peck,kolibiar,cs_cong,halas,reading_order,hl_vop,szk,abbes,hs-char-poly,gk_cong,cf_decomp,njw-hbo}, while the subject of quotients of posets is not covered by many prominent texts on partial orders, such as \cite{dp_itlao,harzheim,stanley_enum,schroeder_os}. General posets thus stand in stark contrast to lattices, where there is a singular natural notion of a congruence, given by preservation of the meet and join operations. Such equivalence relations are known as `lattice congruences' and guarantee that the quotient is itself a lattice, and thus \textit{a fortiori}, a well-defined poset. Indeed, there is a natural notion of congruence for any universal algebra, which specialises to the notion of a lattice congruence, as lattices are universal algebras. Posets, however, are not universal algebras, as they are defined by relations rather than operations.

In defining a class of poset congruences, there is a trade-off between how large the class is, and how much structure the congruences preserve. It is on the one hand desirable to be able to quotient posets by as large a class of equivalence relations as possible, but on the other hand it is desirable for the quotients to preserve a good deal of structure. Naturally, the more equivalence relations one wants to admit, the less structure is in general preserved; and the more structure one wants to preserve, the fewer equivalence relations one can admit.

Let us briefly describe the array of different approaches that exist. One natural approach to take is to allow quotients by as large a class of equivalence relations as possible. In fact, there is a canonical way of constructing a quotient of a poset by an arbitrary equivalence relation \cite{abbes}, notwithstanding what was written earlier. However, in order to make this possible one needs to sacrifice the feature that the elements of the quotient poset should be the equivalence classes of the original poset. Hence, one may wish to restrict to equivalence relations where the elements of the quotient are the equivalence classes of the original poset. The resulting notion was introduced by Sturm \cite{sturm_ao,sturm_vkia,sturm_ecvk,sturm_zao, sturm_cios,sturm_lkim,sturm_loce}, as well as later independently studied in \cite{bj_res,blyth}, and is also related to concepts from \cite{stanley_tpp,trotter_capos}. The downside of this approach is that one is sometimes required to take the transitive closure of the relation one obtains from the quotient. A more well-behaved class of equivalence relations, which do not require the transitive closure to be taken after quotienting is mentioned in \cite{szk} and was independently introduced in \cite{njw-hbo}. But it is natural to desire congruences of posets which preserve yet more structure, such as congruences which preserve upper bounds, as considered in \cite{cs_cong,halas,hl_vop,reading_cambrian,szk}. There are other classes of poset congruences which have particular structure-preserving properties \cite{gk_cong}, such as those that come from lexicographic sums \cite{hh}, or from a group of automorphisms \cite{stanley_quotient_peck,stanley-appl-alg-comb}, those that relate to direct product decompositions of posets \cite{kolibiar}, closures \cite{bj_res,blyth}, or the characteristic polynomial of the poset \cite{hs-char-poly,hallam_applications}. Some types of congruence also relate to specific types of posets \cite{mmp_nat,mmp_dcpo,clp_rel,clp_spp}.

In this survey, we take a pluralist view of congruences and quotients of posets. It is useful to recognise several different types of congruences on posets. Whilst there are some very general notions, many specific examples of congruences possess more structure than is contained in these notions. Indeed, different examples of quotients of posets found in mathematics fall nicely into different specific notions of congruence. Of course, not all types of poset congruence that have appeared so far in the literature may prove equally useful, and there may yet be undiscovered notions of congruence which are very fruitful. To summarise this survey in one sentence: one can quotient any poset by any equivalence relation one wants, but some equivalence relations will preserve more structure than others.

\begin{figure}
\caption{Implications between types of poset congruence. See Section~\ref{sect:comparison} for explanation of terminology}\label{fig:implications}
\[
\begin{tikzpicture}[xscale=3.5,yscale=2.6]

\node(hl) at (-1,0) {\parbox{3cm}{\centering Haviar--Lihov\'a \\ congruence}};
\node(ord) at (-1,-1) {\parbox{2cm}{\centering Order \\ congruence}};
\node(wsc) at (-1,-2) {\parbox{2cm}{\centering $w$-stable \\ congruence}};
\node(gk) at (0.5,-2) {GK congruence};
\node(iii) at (-1,-3) {III-congruence};
\node(clos) at (-0.25,-2.5) {\parbox{2cm}{\centering Closure \\ congruence}};
\node(oac) at (-0.5,-0.5) {\parbox{2cm}{\centering Order-autonomous \\ congruence}};
\node(dcpo) at (-1,-4) {\parbox{3cm}{\centering Natural DCPO \\ congruence}};
\node(woc) at (0,-4) {\parbox{2.5cm}{\centering Weak order \\ congruence}};
\node(cont) at (1,-3.75) {\parbox{2.5cm}{\centering Contraction \\ congruence}};
\node(comp) at (0,-5) {\parbox{2.2cm}{\centering Compatible \\ congruence}};
\node(aer) at (0,-6) {\parbox{4cm}{\centering Arbitrary \\ equivalence relation}};
\node(hom) at (1,-3) {\parbox{2.5cm}{\centering Homogeneous \\ congruence}};
\node(orb) at (1,-0.5) {\parbox{4.5cm}{\centering Orbits of \\ automorphism \\ group}};
\node(kol) at (0.25,-0.5) {\parbox{2cm}{\centering Kolibiar \\ congruence}};

\draw[->,double] (hl) -- (ord) node [midway,left] {\tiny \cite[Thm.~3.5]{hl_vop}};
\draw[->,double] (ord) -- (wsc) node [midway,left] {\parbox{1.8cm}{\tiny \flushright Finite \\ \cite[Thm.~3]{halas}}};
\draw[->,double] (wsc) -- (iii) node [midway,right] {\tiny Strong} node [midway,left] {\tiny \cite[Thm.~4.1]{szk}};
\draw[->,double] (iii) -- (woc);
\draw[->,double] (oac) -- (gk);
\draw[->,double] (woc) -- (comp);
\draw[->,double] (comp) -- (aer);
\draw[->,double] (dcpo) -- (comp);
\draw[->,double] (hom) -- (woc) node [midway,left] {\tiny Finite};
\draw[->,double] (gk) -- (hom) node [midway,right] {\tiny $\{\hat{0}\}$};
\draw[->,double] (ord) -- (gk);
\draw[->,double] (kol) -- (gk) node [midway,left] {\tiny \cite[2.5]{kolibiar}};
\draw[->,double] (orb) -- (gk);
\draw[->,double] (ord) -- (clos) node [midway,left] {\tiny Finite};
\draw[->,double] (clos) -- (woc);
\draw[->,double] (clos) -- (dcpo) node [midway,right] {\tiny \parbox{1.2cm}{\cite[Thm.~2.11]{mmp_nat}}};
\draw[->,double] (cont) -- (comp);
\draw[->,double] (ord) -- (cont) node [midway,right] {\tiny Finite};

\end{tikzpicture}
\]
\end{figure}
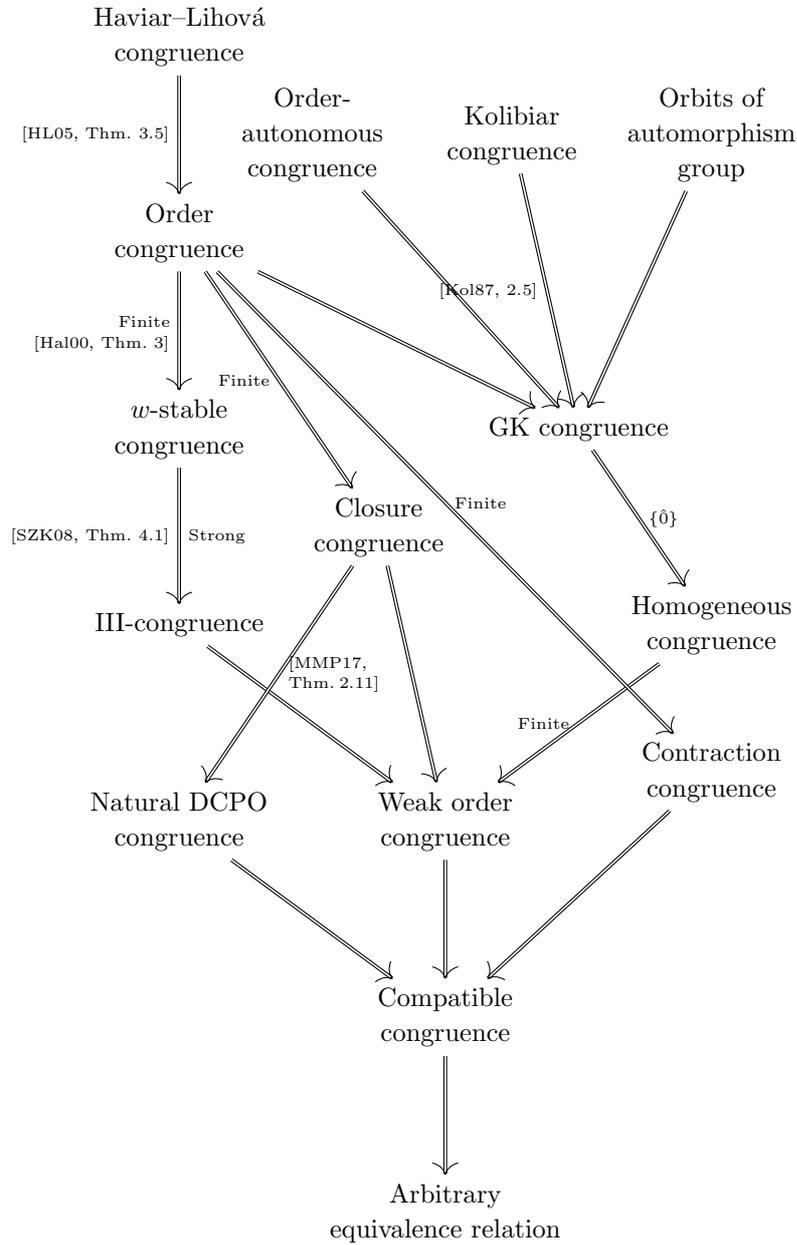

The aims of this paper are thus
\begin{itemize}
\item to describe different approaches taken to quotients of posets in the literature;
\item to give the motivation for the different approaches;
\item to compare and contrast the different approaches: see Figure~\ref{fig:implications} and Table~\ref{table};
\item to provide general philosophical discussion about how one ought to approach taking a quotient of a poset;
\item to survey important examples of quotients of particular posets that have appeared in the literature.
\end{itemize}
The literature on quotients of posets is currently very dispersed, which makes it hard for authors to be aware of what notions exist. Indeed, for many relevant papers, such as \cite{dh_order}, copies are hard to find. Moreover, most papers, understandably, present only their own approach, which can prevent working mathematicians from finding notions of congruences and quotients suited to their own problems.

We begin in Section~\ref{sect:back} by giving background on partially ordered sets. In Section~\ref{sect:quot}, we introduce the problem of taking quotients of posets and discuss possible approaches to this problem in broad terms. We then describe the several different approaches given in the literature in Sections~\ref{sect:not:gen}, \ref{sect:not:lat}, \ref{sect:not:misc}, and~\ref{sect:not:spec}. The types of congruence in Section~\ref{sect:not:gen} are those which aim to admit quotients by many different equivalence relations. In Section~\ref{sect:not:lat}, we examine types of congruence which aim to generalise lattice congruences in some way. Section~\ref{sect:not:misc} then considers congruences which do not fall into these other categories, whilst Section~\ref{sect:not:spec} considers congruences which require additional assumptions on the poset. Following this, in Section~\ref{sect:ex}, we survey important examples of quotients of posets that have appeared in the literature. Finally, in Section~\ref{sect:comparison}, we compare the different approaches.

\afterpage{%
    \clearpage
    \thispagestyle{empty}
    \begin{landscape}
        \centering 
        \renewcommand{\arraystretch}{1.4}%
        \begin{tabular}{l|c|c|c|c|c|c|c|c}
Type of congruence & Self-dual? & \parbox{1.8cm}{\centering Preserves \\ grading?} & \parbox{1.9cm}{\centering Quotient \\ map strong?} & \parbox{1.5cm}{\centering Closed \\ under $\cap$?} & \parbox{1.3cm}{\centering Infinite posets?} & \parbox{2cm}{\centering Lattice \\ congruences?} & \parbox{2cm}{\centering Other \\ requirements} \\
\hline
Equivalence relation & \cmark & \xmark & \xmark & \cmark & \cmark & \xmark & None \\
\hline
Compatible & \cmark & \xmark & \xmark & \cmark\textsuperscript{\eqref{ref:comp:int}} & \cmark & \xmark & None \\
\hline
Weak order & \cmark & \xmark & \cmark & \xmark & \cmark & \xmark & None \\
\hline
III & \xmark & \xmark & \cmark & \xmark & \cmark & \xmark & None \\
\hline
$w$-stable & \cmark & \xmark & \xmark & \cmark\textsuperscript{\eqref{ref:wsc:int}} & \cmark & \cmark & None \\
\hline
Order & \cmark & \xmark & \cmark & \xmark\textsuperscript{\eqref{ref:ord:int}} & \cmark & \cmark & None \\
\hline
Haviar--Lihov\'a & \cmark & \xmark & \cmark & \cmark\textsuperscript{\eqref{ref:hl:int}} & \cmark & \cmark\textsuperscript{\eqref{ref:hl:lat}} & None \\
\hline
GK & \cmark & \xmark & \cmark & \xmark & \cmark & \xmark & None \\
\hline
Order-autonomous & \cmark & \xmark & \cmark & \cmark & \cmark & \xmark & None \\
\hline
Closure & \xmark &\xmark & \cmark & \xmark & \cmark & \xmark & None \\
\hline
Orbits & \cmark & \cmark & \cmark & \xmark & \xmark & \xmark & None \\
\hline
Contraction & \cmark & \xmark & \xmark & \xmark & \cmark & \xmark & None \\
\hline
Kolibiar & \cmark & \xmark & \cmark\textsuperscript{\eqref{ref:kol:full}} & \xmark & \cmark & \cmark\textsuperscript{\eqref{ref:kol:lat}} & Directed \\
\hline
Homogeneous & \xmark & \xmark & \cmark & \xmark & \xmark & \xmark & Unique min. \\
\hline
Natural DCPO & \xmark & \xmark & \xmark & \xmark & \cmark & \xmark & DCPO \\
\hline
		\end{tabular}
		\bigskip
        \captionof{table}{Comparing different poset congruences. See Section~\ref{sect:comparison} for explanation of terminology}\label{table}
        \renewcommand{\arraystretch}{1}%
    \end{landscape}
    \clearpage
}

\section{Background}\label{sect:back}

We begin by giving basic definitions for partial orders. 

\subsection{Partially ordered sets}\label{sect:back:poset}

Given a set $\set$, a \emph{relation} $\rel$ on $\set$ is a subset of the Cartesian product $\set \times \set$. We will also sometimes write that $(\set, \rel)$ is a relation, if we want to make both the symbol for the relation and the symbol for the underlying set clear. Somewhat confusingly, it is usual also to refer to the elements of $R$ as relations. If $(x, y) \in \rel$, then we write $x \rel y$. A relation $\rel$ is
\begin{itemize}
\item \emph{reflexive} if $x \rel x$ for all $x \in \set$, 
\item \emph{symmetric} if $x \rel y$ implies that $y \rel x$,
\item \emph{anti-symmetric} if $x \rel y$ and $y \rel x$ together imply that $x = y$,
\item \emph{transitive} if $x \rel y$ and $y \rel z$ together imply that $x \rel z$.
\end{itemize}
A relation $\rel$ on a set $\set$ is called a \emph{partial order} if it is reflexive, anti-symmetric, and transitive. Here we call $\set$ a \emph{partially ordered set}, or \emph{poset}. We usually write $P$ instead of $\set$ and $\leqslant$ instead of $\rel$ if we have a partially ordered set, rather than only a set with a relation, so that we write that $(P, \leqslant)$ is a poset. We will often simply write that $P$ is a poset, meaning that there is a partial order on $P$ denoted~$\leqslant$. As is very standard, given $p, q \in P$, we write
\begin{itemize}
\item $p < q$ if $p \leqslant q$ and $p \neq q$,
\item $p \not\leqslant q$ (resp.\ $p \not< q$) if $p \leqslant q$ (resp.\ $p < q$) is not the case, and
\item $p \geqslant q$ (resp.\ $p > q$) if $q \leqslant p$ (resp.\ $q < p$).
\end{itemize}

\newpage

The \emph{transitive closure} $\trans{R}$ of a relation $R$ on $S$ is the smallest transitive relation containing $R$ or, equivalently, the intersection of all of the transitive relations containing $R$. Reflexive and symmetric closures are defined likewise. Given $x, z \in P$ such that $x < z$, we say that $z$ \emph{covers} $x$ if there is no $y \in P$ such that $x < y < z$. In this case we write $x \lessdot z$ and refer to this as a \emph{covering relation}. A finite poset is equal to the transitive closure of its covering relations, but this is not always true for infinite posets.

Given a subset $A \subseteq P$, an \emph{upper bound} of $A$ is an element $u \in P$ such that $a \leqslant u$ whenever $a \in A$. The notion of a \emph{lower bound} is defined dually. We denote by
\begin{align*}
\ls{P}{A} &= \spt{ p \in P \st p \leqslant a \text{ for all } a \in A}, \\
\us{P}{A} &= \spt{ p \in P \st p \geqslant a \text{ for all } a \in A}
\end{align*}
the respective sets of lower and upper bounds of $A$ in $P$. Given a subset $A = \{a_{1}, a_{2}, \dots, a_{k}\}$, we may sometimes write $L_{P}(a_{1}, a_{2}, \dots, a_{k})$ and $U_{P}(a_{1}, a_{2}, \dots, a_{k})$ for $L_{P}(A)$ and $U_{P}(A)$. A \emph{directed poset} is a poset where any finite subset has an upper bound. The \emph{supremum} or \emph{least upper bound} of $A$, if it exists, is defined to be an element $s$ such that for any upper bound $u$ of $A$, we have that $s \leqslant u$. It is easy to see that, if $A$ has a supremum $s$, then $s$ is unique. We denote the supremum of $A$ by $\sup A$, if it exists. The \emph{infimum} or \emph{greatest lower bound} $\inf A$ is defined dually. Given a poset $P$, a \emph{minimal element} of $P$ is an element $m$ such that $p \not< m$ for all $p \in P$. A \emph{maximal element} of $P$ is defined dually.

A subset $I$ of $P$ is an \emph{interval} if it is of the form $I = [p,r] := \spt{q \in P \st p \leqslant q \leqslant r}$. A subset $P' \subset P$ is called \emph{convex} if whenever $p, r \in P'$ and $p \leqslant q \leqslant r$, then $q \in P'$.

A poset $P$ is a \emph{totally ordered set} and $\leqslant$ is a \emph{total order} if for any $p, q \in P$, we either have $p \leqslant q$ or $p \geqslant q$. That is, in a totally ordered set, any two elements are \emph{comparable}. A \emph{chain} $C$ in $P$ is a subset totally ordered by $\leqslant$. An \emph{antichain} of $P$ is a subset of $P$ such that no two distinct elements are comparable.

Given a poset $(P, \leqslant)$, the \emph{dual poset} is the poset $(\overline{P}, \overline{\leqslant})$ where $\overline{P} = P$ and $p \leqslant q$ if and only if $q \mathrel{\overline{\leqslant}} p$.

\subsubsection{Maps between posets}

Given posets $P$ and $Q$, a map $f \colon P \to Q$ is said to be \emph{order-preserving} if $f(p) \leqslant f(q)$ whenever $p \leqslant q$. The map $f$ is said to be \emph{order-reversing} if $f(p) \geqslant f(q)$ whenever $p \leqslant q$. Order-preserving maps are also called \emph{isotone} and order-reversing maps are also called \emph{antitone}.

An \emph{isomorphism} of posets is an order-preserving bijection whose inverse is order-preserving; an \emph{anti-isomorphism} of posets is an order-reversing bijection whose inverse is order-reversing.
If there is an isomorphism of posets $P \to Q$, we write $P \cong Q$.
A poset is \emph{self-dual} if it is isomorphic to its dual. An \emph{automorphism} of $P$ is an isomorphism $P \to P$; an \emph{anti-automorphism} is an anti-isomorphism $P \to P$. Anti-isomorphisms and anti-automorphisms are also called \emph{dual isomorphisms} and \emph{dual automorphisms}, respectively.

A map $f\colon P \to Q$ is \emph{strong} if it is order-preserving and if whenever $f(p) \leqslant f(p')$ there exist $\hat{p},\hat{p}' \in P$ such that $\hat{p} \leqslant \hat{p}'$ and $f(\hat{p}) = f(p)$ and $f(\hat{p}') = f(p')$. Strong maps can be thought of as being surjective on the relations. A poset $P$ is a \emph{subposet} of a poset $Q$ if there is a strong injection $P \hookrightarrow Q$.

\subsubsection{Graded posets}\label{sect:back:def:graded}

A poset $P$ is \emph{graded} if there is a \emph{rank function} $\rho\colon P \to \mathbb{Z}_{>0}$ such that if $p < q$ in $P$ then $\rho(p) < \rho(q)$ and if $p \lessdot q$, then $\rho(q) = \rho(p) + 1$. The value of $\rho(p)$ is called the \emph{rank} of $p$. We say that $P$ has \emph{rank} $n$ if the largest value of $\rho(p)$ is $n$ and the lowest value of $\rho(p)$ is $0$. The \emph{ranks} of $P$ are the subsets $P_{i} = \spt{p \in P \st \rho(p) = i}$.

\subsection{Equivalence relations}

A relation $\rel$ on a set $\set$ is called an \emph{equivalence relation} if it is reflexive, symmetric, and transitive. In this paper, we usually denote equivalence relations by $\equ$. A \emph{partition} of $\set$, is a set $\{\set_{i} \st i \in I\}$ of pairwise disjoint non-empty subsets of $\set$ such that $\set = \bigcup_{i \in I} \set_{i}$, where $I$ is some indexing set. Partitions of $\set$ are equivalent to equivalence relations on $\set$. Namely, given an equivalence relation $\equ$ on a set $\set$, there is a partition of $\set$ into sets $\set_{i}$, where, if $s \in \set_{i}$, then $t \in \set_{i}$ if and only if $s \equ t$. In this case, $\set_{i}$ is called the \emph{equivalence class} of $s$ and is denoted $[s]$ or $[s]_{\equ}$. Hence, $t \in [s]$ if and only if $t \equ s$. Conversely, given a partition $\spt{\set_{i} \st i \in I}$ of $\set$, where we write $[s] = \set_{i}$ for the unique $\set_{i}$ such that $s \in \set_{i}$, the corresponding equivalence relation is given by $s \equ t$ if and only if $[s] = [t]$. The set of $\equ$-equivalence classes of $S$ is denoted by $S/{\equ} = \{S_{i}\}_{i \in I}$.

\subsection{Pre-orders}\label{sect:back:preord}

A relation $R$ on a set $S$ which is reflexive and transitive is known as a \emph{pre-order}. Pre-orders are thus a simultaneous generalisation of partial orders and equivalence relations, with the former being anti-symmetric pre-orders and the latter being symmetric pre-orders.

There is a canonical way of defining a partial order from a pre-order. Indeed, let $(P, \pleq)$ be a pre-order. There is an equivalence relation $\equ$ on $P$ given by $p \equ q$ if and only if $p \pleq q$ and $p \pgeq q$. Then one defines the \emph{collapse} of the pre-order $(P, \pleq)$ to be the poset $\collapse{P}$ with underlying set the set of $\equ$-equivalence classes $\quot$ and relation given by $[p] \leqslant [q]$ if and only if $p' \pleq q'$ for some $p' \in [p]$ and $q' \in [q]$. The following result is well-known and straightforward.

\begin{proposition}
The collapse $\collapse{P}$ is a well-defined poset. Moreover, for $[p], [q] \in \collapse{P}$, we have $[p] \leqslant [q]$ if and only if $p' \pleq q'$ for all $p' \in [p]$ and $q' \in [q]$.
\end{proposition}

\begin{remark}
The collapse operation gives a functor from the category $\preordcat$ of pre-orders to the category $\posetcat$ of posets, where the morphisms in each category are given by order-preserving maps. In fact, it is the left adjoint in an adjunction with the forgetful functor $\mathsf{forget}$ from $\posetcat$ to $\preordcat$. This means that for any pre-order $P$ and poset $Q$, there is a bijection \[\Hom_{\posetcat}(\collapse{P}, Q) \cong \Hom_{\preordcat}(P, \mathsf{forget}(Q))\] which is natural in both $P$ and $Q$. For more details on adjunctions, see \cite[Chapter~2]{leinster}.
\end{remark}

\subsection{Lattices}

A lattice is a poset with certain properties which give rise to additional algebraic structure. According to Birkhoff \cite[p.\ 6]{birkhoff}, the concept of a lattice was first studied in depth by Dedekind, under the name `Dualgruppe' \cite[pp.\ 113--4]{dedekind_lat}, see also \cite{dd_vuz,dedekind_drei}. Partial orders then arose from lattices, an early reference to partial orders being Hausdorff's book on set theory \cite[Sechstes~Kapitel, \S1]{haus_gdm}.

A \emph{lattice} is a partially ordered set $L$ such that for every pair of elements $x, y \in L$, $\{x, y\}$ has both a supremum and an infimum. Here the supremum is called the \emph{join} and is denoted $x \join y$, and the infimum  is called the \emph{meet} and is denoted $x \meet y$. A \emph{complete lattice} is a poset~$L$ such that any subset $X \subseteq L$ has a infimum $\bigwedge X$ and a supremum $\bigvee X$. A poset is a \emph{join-semilattice} if every pair of elements has a join; dually, a poset is a \emph{meet-semilattice} if every pair of elements has a meet.

Note that meet and join are binary operations on a lattice. These lattice operations actually suffice to define the partial order. Indeed, we have the following theorem.

\begin{theorem}[{\cite[Theorem 8]{birkhoff}}]\label{thm:lattice_operations}
If $L$ is a set with two binary operations $\meet$ and $\join$ such that
\begin{enumerate}[label=\textup{(}\arabic*\textup{)}]
\item $x \meet x = x$, \quad $x \join x = x$;\label{op:lattice:refl}
\item $x \meet y = y \meet x$, \quad $x \join y = y \join x$;\label{op:lattice:sym}
\item $x \meet (y \meet z) = (x \meet y) \meet z$, \quad $x \join (y \join z) = (x \join y) \join z$;\label{op:lattice:ass}
\item $x \meet (x \join y) = x \join (x \meet y) = x$;\label{op:lattice:swap}
\end{enumerate}
then $L$ is a lattice with $\meet$ and $\join$ the meet and join operations.The partial order on $L$ is defined by $x \leqslant y$ if and only if $x \meet y = x$ or, alternatively, $x \join y = y$.
\end{theorem}

\subsubsection{The Dedekind--MacNeille completion}

The Dedekind--MacNeille completion of a poset will be important at various points in this survey. This is the smallest complete lattice which contains a given poset. It was introduced in \cite{macneille}, inspired by, and generalising, Dedekind's construction of the real numbers from the rationals \cite{dedekind_cut}.

The \emph{Dedekind--MacNeille completion} of a poset $P$ has underlying set \[\mac{P} = \spt{A \subseteq P \st \ls{P}{\us{P}{A}} = A},\] with order relation given by inclusion. The natural embedding \[P \hookrightarrow \mac{P}\] is given by $p \mapsto \ls{P}{p}$.

\section{Quotients of posets}\label{sect:quot}

In this section we outline the general problem with taking the quotient of a poset by an equivalence relation, and consider how one may deal with this problem. We contrast this with the situation of lattices.

\subsection{The essential problem}\label{sect:quot:prob}

In general, if $S$ is some mathematical structure, then the quotient of $S$ by an equivalence relation $\equ$ has the set of $\equ$-equivalence classes $S/{\equ}$ as its underlying set. If the equivalence relation is well-behaved, then the set $S/{\equ}$ will inherit the mathematical structure from $S$ in a natural way.

Given a poset $(P, \leqslant)$ and an equivalence relation $\equ$, the most natural way to define the quotient poset $(\quot, \qrel)$ is by defining $\qrel$ such that, given $[p], [q] \in \quot$, we have that $[p] \qrel [q]$ if and only if there exists $p' \in [p]$ and $q' \in [q]$ such that $p' \leqslant q'$. One way of obtaining this definition of the quotient relation is to note that there is a natural map of sets $P \to \quot$ given by $p \mapsto [p]$. We refer to this as the \emph{quotient map}. Requiring that the quotient map be order-preserving then produces the definition of the quotient relation $\qrel$ given.

The quotient relation $\qrel$ constructed in this way is in general only a reflexive relation on $\quot$. It is not generally transitive  or anti-symmetric, as can be seen from the following examples. This is the key point which demands that we seek a class of equivalence relations that give well-defined quotient posets, or take a different approach.

\begin{example}\label{ex:not_antisym}
We give an example to show that the relation $\qrel$ is not generally anti-symmetric. Consider the poset $P = \{p < q < r\}$ and the equivalence relation $\equ$ given by the partition $\{p,r\}, \{q\}$. Here the quotient relation $\qrel$ on $\quot$ is equal to the reflexive closure of the relation $\{([p], [q]), ([q], [r])\}$, which is not a partial order since $[p] \qrel [q] \qrel [r]$, when $[p] = [r] \neq [q]$.
\end{example}

\begin{example}\label{ex:not_trans}
We give an example to show that the relation $\qrel$ is not generally transitive. Consider the partial order $\leqslant$ on the set $P = \{p, q, q', r\}$ given by the reflexive closure of $\{(p, q), (q', r)\}$. Let $\equ$ be the equivalence relation given by the partition $\{p\}, \{q, q'\}, \{r\}$. Then the quotient relation $\qrel$ on $\quot$ is equal to the reflexive closure of the relation $\{([p], [q]), ([q'], [r])\}$. However, this is not a partial order, because $[p] \qrel [q]$ and $[q] = [q'] \qrel [r]$, but $[p] \not\qrel [r]$.
\end{example}

Hence, it is desirable to restrict the equivalence relations $\equ$ which can be used to take a quotient of the poset $P$. One may then call the members of this restricted class of equivalence relations `congruences'. In the literature several different restrictions on the equivalence relations are imposed. Another approach, of course, is to give a different definition of the quotient of a poset by a particular equivalence relation. This approach is taken in Section~\ref{sect:not:univ}.

\begin{remark}\label{rmk:graphs}
The definition of the quotient also applies to arbitrary relations and directed graphs, as in \cite{quotient_digraphs}.
\end{remark}

\subsubsection{Morphism perspective}

Given an order-preserving map of posets $f\colon P \to Q$, the \emph{kernel} of $f$ is the equivalence relation $\equ$ where $p \theta p'$ if and only if $f(p) = f(p')$. Another perspective on congruences of posets is that all congruences on posets should be kernels of some order-preserving map, namely the map $P \to \quot$. Instead of characterising a class of equivalence relations directly, one can look for a class of maps. One can then take congruences to be the kernels of this particular class of maps. We often refer to this as the `morphism perspective'.

\subsection{Lattice congruences}

The situation for considering quotients of lattices is altogether better than that of quotients of posets. This is because lattices possess not only the order \emph{relation}, but also the meet and join \emph{operations}. The natural class of well-behaved equivalence relations are those that respect the meet and join operations. Such equivalence relations produce well-defined quotient lattices. More detail on lattice congruences can be found in \cite{birkhoff,dp_itlao,roman,graetzer}.

Indeed, given a lattice $L$, a \emph{lattice congruence} on $L$ is an equivalence relation $\equ$ such that whenever $x \equ y$ in $L$, we have $(x \join z) \equ (y \join z)$ and $(x \meet z) \equ (y \meet z)$ for all $z \in L$. If $L$ is a complete lattice, then a \emph{complete lattice congruence} on $L$ is an equivalence relation $\equ$ such that for any indexing set $I$ and subsets $\{x_{i}\}_{i \in I}, \{y_{i}\}_{i \in I} \subseteq L$ with $x_{i} \equ y_{i}$ for all $i$, then \[\bigwedge\{x_{i} \st i \in I\} \mathrel{\equ} \bigwedge\{y_{i} \st i \in I\}\] and \[\bigvee\{x_{i} \st i \in I\} \mathrel{\equ} \bigvee\{y_{i} \st i \in I\}.\]
If these properties are only satisfied with respect to one of meet or join, then the equivalence relation is called a \emph{\textup{(}complete\textup{)} meet-semilattice conguence} or a \emph{\textup{(}complete\textup{)} join-semilattice conguence} respectively.

\begin{remark}
Lattices are `universal algebras', since they are defined by a set of elements and operations. The definition of lattice congruence corresponds to the general definition of congruence for universal algebras, which states that operations applied to equivalent elements must give the same result. For more details, see \cite[Chapter~VI]{birkhoff}.

Posets, on the other hand, are not universal algebras, since they are not defined by operations, but rather by a relation. Nevertheless, one could try to adapt the notion of a congruence for universal algebras to posets, by requiring that relations continue to hold if one substitutes an element for an equivalent one; that is, if $p \leqslant q$ and $p \equ p'$, then $p' \leqslant q$, and that if $p \leqslant q$ and $q \equ q'$, then $p \leqslant q'$.
One can check that quotients of posets by such congruences give well-defined posets.
In fact, the resulting notion is that of an `order-autonomous congruence' discussed in Section~\ref{sect:not:oac}.
However, this notion of congruence is much too restrictive to capture all interesting cases of quotient posets in mathematics.
\end{remark}

The quotient of a lattice by a lattice congruence is always a lattice, and hence a partial order.

\begin{proposition}[{\cite{birkhoff}}]
If $\equ$ is a lattice congruence on a lattice $L$, then the quotient $\lquot$ is a lattice and, \emph{a fortiori}, a partial order.
\end{proposition}

One can then show that the partial order $\qrel$ on $\lquot$ for a lattice $L$ given by applying Theorem~\ref{thm:lattice_operations} to $\lquot$ is the same as would be given by considering $\lquot$ as the quotient of the poset~$L$.

Usually, the easiest way to verify whether a given equivalence relation on a lattice is a lattice congruence is to apply the following lemma.

\begin{lemma}[{\cite[Lemma~4]{gs_icr}}]\label{lem:tech}
An equivalence relation $\equ$ on a lattice $L$ is a congruence relation if and only if the following three properties are satisfied for all $w, x, y, z \in L$.
\begin{enumerate}[label=\textup{(}\arabic*\textup{)}]
\item $x \equ y$ if and only if $(x \meet y) \equ (x \join y)$.\label{op:tech:collapse}
\item If $x \leqslant y$ and $x \equ y$, then $(x \meet w) \equ (y \meet w)$ and $(x \join w) \equ (y \join w)$.\label{op:tech:subs}
\end{enumerate}
\end{lemma}

Note that there are several other criteria for an equivalence relation to be a lattice congruence \cite{dorfer,reading_order}.

\subsection{Discussion}

We seek a class of equivalence relations $\equ$ on posets $P$ which produce well-defined quotients. Let us briefly consider the desiderata for such a class of equivalence relations.

Clearly, the principal requirement is that the quotient relation $\qrel$ is a partial order. One approach would be to take this as the only requirement, thereby obtaining the largest class of equivalence relations producing well-defined quotient posets. This is the approach taken in Section~\ref{sect:not:njw}. One can, in fact, take an even weaker approach than this by requiring only that the transitive closure of the quotient relation is a partial order, as in Section~\ref{sect:not:comp}.

However, there are many reasons to study stronger notions of poset congruences. It is reasonable to be interested in equivalence relations on posets which preserve more than just the partial order structure, but which may also preserve upper and lower bounds, for instance. Such approaches are considered in Section~\ref{sect:not:lat}. The point is that congruences of posets that arise in natural mathematical examples often have stronger properties than simply producing a well-defined quotient, so it is important to have types of congruence which take account of these stronger properties. Such types of congruence are also easier to work with.

A related desire concerns how a particular class of equivalence relations interacts with lattices and lattice congruences. The idea is that lattice congruences are clearly the right notion of congruence for lattices, so one ought to try to extend the notion of a lattice congruence to all posets. Hence, one might require that the class of poset congruences considered coincides with the class of lattice congruences when restricted to lattices. This is the case for most of the notions of congruence studied in Section~\ref{sect:not:lat}. However, one might feel that, for a sufficiently general class of poset congruences, one ought to be able to produce quotients of lattices which are not themselves lattices. Indeed, there are interesting examples of this given in Section~\ref{sect:ex}.

There are many other specific circumstances in which interesting classes of congruences arise, such as those which result from closure operators or groups of automorphisms. These sorts of congruences are considered in Section~\ref{sect:not:misc}. Similarly, specific subclasses of posets may possess classes of congruences which are particular to them, as considered in Section~\ref{sect:not:spec}. Indeed, lattices give such an example.

We end with a brief list of the different desiderata used to motivate various types of poset congruence in the literature. It is impossible to simultaneously satisfy all of these desiderata.
\begin{itemize}
\item The notion includes as many equivalence relations as possible.
\item Constructing the quotient poset does not require taking the transitive closure.
\item The notion coincides with that of a lattice congruence in the case of lattices.
\item The quotients preserve upper bounds, or other order-theoretic features of the poset.
\item The definition is natural and not too complicated.
\item The definition is relatively easy to check in practice.
\item There exist examples occurring in nature.
\item The intersection of two congruences should be a congruence.
\item The set of all congruences on a poset should itself possess some nice structure, such as being a lattice or a complete lattice.
\end{itemize}

We now survey the different notions of congruences of posets that have appeared in the literature, grouping similar notions of congruence together.

\section{General notions}\label{sect:not:gen}

We begin by considering notions of poset congruence from the literature which aim to incorporate as large a class of equivalence relations as possible, perhaps subject to some minimal restrictions.

\subsection{Universal property approach}\label{sect:not:univ}

The most general approach to taking quotients of posets allows quotients with respect to any equivalence relation by defining the quotient using a universal property. This approach can be found in \cite{abbes}.

\begin{definition}\label{def:univ}
Given a poset $P$ and an equivalence relation $\equ$, a \emph{universal quotient of $P$} by $\equ$ is any pair $(Q, f)$, where $Q$ is a poset and $f \colon P \to Q$ an order preserving map, which satisfies the following universal property. If $R$ is a poset and $g \colon P \to R$ is an order-preserving map which is constant on the $\equ$-equivalence classes of $P$, then there is a unique order-preserving map $h \colon Q \to R$ such that $g = h \circ f$.
\end{definition}

One can construct a universal quotient for any equivalence relation, and such a universal quotient is unique up to unique isomorphism.

\begin{proposition}[\cite{abbes}]\label{prop:univ:prop}
Let $P$ be a poset with $\equ$ an equivalence relation on~$P$.
\begin{enumerate}[label=\textup{(}\arabic*\textup{)}]
\item There is a poset $Q$ with an order-preserving map $f \colon P \to Q$ such that $(Q, f)$ is the universal quotient of $P$ with respect to $\equ$.\label{op:univ:prop:exists}
\item Given another universal quotient $(Q', f')$ of $P$ by $\equ$, there is a unique isomorphism $h \colon Q \to Q'$ such that $f' = h \circ f$.\label{op:univ:prop:unique}
\end{enumerate}
\end{proposition}

Here \ref{op:univ:prop:unique} is proven by a standard argument, whilst \ref{op:univ:prop:exists} uses the following construction.

\begin{construction}[\cite{abbes}]
We construct the universal quotient of a poset $P$ by an equivalence relation $\equ$. We first let $\quot$ be the set of $\equ$-equivalence classes of $P$, as ever. There is an obvious map of sets $f'\colon P \to \quot$. We define a pre-order $\pleq$ on $\quot$ by specifying that if $p \leqslant q$ in $P$ then $[p] \pleq [q]$ in $\quot$, and taking the transitive closure $\trans{\pleq}$ of the resulting relation. We then define $Q = \collapse{\quot, \trans{\pleq}}$, the collapse of the pre-order $(\quot, \trans{\pleq})$ from Section~\ref{sect:back:preord}. The map $f \colon P \to Q$ is defined to be $f'' \circ f'$ where $f''$ is the canonical map $f'' \colon \quot \to Q$. We then have that $f \colon P \to Q$ is the universal quotient of $P$ by $\equ$.
\end{construction}

\begin{example}\label{ex:univ_prop}
We return to the example given in Example~\ref{ex:not_antisym}, where we have the poset $P = \{p < q < r\}$ under the equivalence relation given by the partition $\{p, r\}, \{q\}$. The pre-order on $\quot$ is given by $\{p, r\} \pleq \{q\} \pleq \{p, r\}$, and the collapse of this pre-order is the one-element poset $Q = \{P\}$.
\end{example}

\subsection{Compatible congruences}\label{sect:not:comp}

The principal drawback of the universal property approach to quotients of posets is that the kernel of the map $f \colon P \to Q$ does not always coincide with the original equivalence relation $\equ$, as in Example~\ref{ex:univ_prop}. Hence, really what we are doing when we take the universal quotient by $\equ$ is replacing the original equivalence relation $\theta$ with some better behaved equivalence relation $\widehat{\equ}$. It makes sense to characterise these better behaved equivalence relations. They are given by the following notion, first considered by Sturm \cite{sturm_ao,sturm_vkia,sturm_ecvk,sturm_zao, sturm_cios,sturm_lkim,sturm_loce}. We use the definition from \cite{krs}, which is in turn inspired by \cite{cl_cong}. Here, and henceforth in this paper, we use the notation $[n] = \{1, 2, \dots, n\}$.

\begin{definition}[{\cite[Definition~2.1]{krs}}]\label{def:compatible}
Let $P$ be a poset with $\equ$ an equivalence relation.
\begin{enumerate}
\item A sequence $(p_{0}, p_{1}, \dots, p_{n})$ is called a \emph{$\equ$-sequence} if for each $i \in [n]$, either $p_{i - 1} \equ p_{i}$ or $p_{i - 1} < p_{i}$. If, additionally, $p_{0} = p_{n}$, then $(p_{0}, p_{1}, \dots, p_{n})$ is called a \emph{$\equ$-circle}.
\item We say that $\equ$ is a \emph{compatible congruence} if we have that $[p_{0}] = [p_{1}] = \dots = [p_{n}]$ for any $\equ$-circle $(p_{0}, p_{1}, \dots, p_{n})$.
\end{enumerate}
\end{definition}

\begin{remark}
This notion was also considered in \cite{bj_res,blyth} and coincides with Trotter's notion of an order-preserving partition, which was introduced in a different context \cite{trotter_capos}. It is also the same as the notion of `compatibility' used by Stanley in the context of order polytopes \cite{geissinger,stanley_tpp}, hence the name we choose.

The papers \cite{cl_cong,cl_class} consider compatible congruences in the context of partially ordered universal algebras. The relevant notion of a congruence on such algebras is that of a congruence of the universal algebra which is compatible with the partial order in the sense of Definition~\ref{def:compatible}.
\end{remark}

Compatible congruences have the following properties.

\begin{proposition}[{\cite[Theorem~3.2]{krs},\cite{sturm_vkia,sturm_lkim,cl_cong}}]\label{prop:comp:prop}
Let $P$ be a poset.
\begin{enumerate}[label=\textup{(}\arabic*\textup{)}]
\item An equivalence relation $\equ$ on $P$ is compatible if and only if any of the following hold:
	\begin{enumerate}[label=\textup{(}\alph*\textup{)}]
	\item the transitive closure $\trans{\qrel}$ of the quotient relation $\qrel$ on $\quot$ is anti-symmetric and hence a partial order;
	\item $\leqslant$ can be extended to a total order $\leqslant_{t}$ on $P$ such that the $\equ$-equivalence classes are intervals with respect to $\leqslant_{t}$; or
	\item $\leqslant$ can be extended to a pre-order $\rho$ on $P$ such that $\equ = \rho \cap \overline{\rho}$, where $\overline{\rho}$ denotes the opposite pre-order.
	\end{enumerate}
\item If $\spt{\equ_{i} \st i \in I}$ is a set of compatible congruences on $P$, then $\equ = \bigcap_{i \in I} \equ_{i}$ is a compatible congruence on~$P$. Indeed, the set of compatible congruences on $P$ forms a complete lattice, where the meet is given by intersection.\label{op:comp:comp_lat}
\item If $\equ$ is a compatible congruence on $P$, then the $\equ$-equivalence classes are convex sets.
\end{enumerate}
\end{proposition}

All of the types of congruence in the subsequent sections will also be compatible, and so will have convex equivalence classes too. Note that having convex equivalence classes is not in general sufficient for an equivalence relation to be a compatible congruence \cite[Fig.~1]{stanley_tpp}.

\begin{remark}
The complete lattice of compatible congruences does not in general form a sublattice of the lattice of equivalence relations, despite Proposition~\ref{prop:comp:prop}\ref{op:comp:comp_lat}, since the join of two compatible congruences is not usually simply their join as equivalence relations, but rather the smallest compatible congruence containing this.
\end{remark}

Compatible congruences are related to universal quotients as follows. Note that Proposition~\ref{prop:comp:prop}\ref{op:comp:comp_lat} gives a well-defined smallest compatible congruence~$\widehat{\theta}$ containing a given equivalence relation~$\theta$, obtained by taking the intersection of all compatible congruences containing~$\theta$. There is always at least one compatible congruence containing~$\theta$, namely, the equivalence relation with one equivalence class. Of course, when we talk about containment of equivalence relations on $P$, we mean containment as subsets of $P \times P$.

\begin{proposition}\label{prop:univ->comp}
Let $P$ be a poset.
\begin{enumerate}[label=\textup{(}\arabic*\textup{)}]
\item Let $\equ$ be an equivalence relation on $P$. If $f \colon P \to Q$ is the universal quotient of $P$ by $\equ$, then $Q \cong (P/\widehat{\equ}, \trans{\leqslant_{\widehat{\equ}}})$, the quotient of $P$ by $\widehat{\equ}$, where $\widehat{\equ}$ is the smallest compatible congruence containing $\equ$.
\item Consequently, if $\equ$ is a compatible congruence on $P$ and $f \colon P \to Q$ is the universal quotient of $P$ by $\equ$, then $\equ = \ker f$.
\end{enumerate}
\end{proposition}

Proposition~\ref{prop:univ->comp} establishes that one can restrict quotients of posets to quotients of posets by compatible congruences, since quotienting by an arbitrary equivalence relation is equivalent to quotienting by the smallest compatible congruence containing it.

Recall that one may also consider quotients of posets from the `morphism perspective' discussed in Section~\ref{sect:quot}. Here, one wants to characterise the set of maps which have the set of congruences as their kernels. The morphism perspective on compatible congruences simply gives them as the kernels of order-preserving maps.

\begin{proposition}[{\cite{sturm_lkim}}]
An equivalence relation $\equ$ on a poset $P$ is a compatible congruence if and only if it is the kernel of an order-preserving map $f \colon P \to Q$.
\end{proposition}

In the literature, there has been interest in quotienting graded posets in such a way that they remain graded, for instance \cite[Lemma~7]{hs-char-poly}. This can be ensured by requiring that the rank function is constant on equivalence classes. All congruences for which this is true are compatible, and the quotient poset is graded in the natural way.

\begin{proposition}\label{prop:grade&comp}
Let $P$ be a graded poset with rank function $\rho$. Let further $\equ$ be an equivalence relation on $P$ such that $\rho(p) = \rho(q)$ whenever $p \equ q$.
\begin{enumerate}[label=\textup{(}\arabic*\textup{)}]
\item The equivalence relation $\equ$ is a compatible congruence.
\item Moreover, $\quot$ is a graded poset with rank function $\tilde{\rho}$ where $\tilde{\rho}([p]) = \rho(p)$ for all $p \in P$.
\end{enumerate}
\end{proposition}

If one had to keep only one type of congruence on posets, it would be compatible congruences, since by Proposition~\ref{prop:univ->comp}, they provide the largest possible class of quotients of posets. However, the properties of compatible congruences are still quite weak, and, as we shall see, it is natural to consider congruences which preserve more structure.

\begin{remark}
In \cite{js}, the authors study the order complex of the lattice of compatible congruences on a poset with $n$ elements, where $n \geqslant 3$. They find that this complex is homotopy equivalent to a wedge of $(n-3)$-spheres and that, if $P$ is connected, then the number of spheres is equal to the number of linear extensions of $P$. A consequent of this is that In \cite{sw}, it is shown that the lattice of compatible congruences on a poset $P$ is always CL-shellable. Each of these results therefore entails that the lattice of compatible congruences is Cohen--Macaulay \cite{baclawski,bgs_cm}.
\end{remark}

\subsection{Weak order congruences}\label{sect:not:njw}

The shortcoming of compatible congruences is that, in general, one needs to take the transitive closure of the quotient relation in order to obtain a well-defined quotient poset. It is natural to consider the more well-behaved class of congruences for which this is not necessary. This notion was first noted in \cite{szk} under the name `II-congruences', and was later independently rediscovered in \cite{njw-hbo}.

\begin{definition}[\cite{njw-hbo}]\label{def:njw}
Let $P$ be a poset with $\equ$ an equivalence relation on $P$. Then we say that $\equ$ is a \emph{weak order congruence} if
\begin{enumerate}
\item given $p, q, p', q' \in P$ such that $p \leqslant q$, $q \equ q'$, $q' \leqslant p'$, and $p' \equ p$, then $p \equ q$, and\label{op:njw:anti-sym}
\item given $p, q, q', r \in P$ such that $p \leqslant q$, $q \equ q'$, and $q' \leqslant r$, then there exist $p', r' \in P$ such that $p \equ p'$, $p' \leqslant r'$, and $r' \equ r$.\label{op:njw:trans}
\end{enumerate}
\end{definition}

The defining conditions of a weak order congruence $\equ$ thus simply amount to specifying that the quotient relation $\qrel$ is a partial order. \eqref{op:njw:anti-sym} holds if and only if $\qrel$ is anti-symmetric and \eqref{op:njw:trans} holds if and only if $\qrel$ is transitive.

\begin{proposition}[{\cite{szk,njw-hbo}}]\label{prop:njw:props}
\begin{enumerate}[label=\textup{(}\arabic*\textup{)}]
\item If $\equ$ is an equivalence relation on a poset $P$, then the quotient relation $\qrel$ on $\quot$ is a partial order if and only if $\equ$ is a weak order congruence.\label{op:njw:works}
\item An equivalence relation on a poset $P$ is a weak order congruence if and only if it is the kernel of a strong map $f \colon P \to Q$.\label{op:njw:morph_pers}
\end{enumerate}
\end{proposition}

Part~\ref{op:njw:morph_pers} gives the morphism perspective on weak order congruences. One way of viewing this result is that weak order congruences correspond to maps which are surjective both on the level of objects and on the level of relations. This is an intuitive notion of a quotient of a poset.

\section{Notions inspired by lattice congruences}\label{sect:not:lat}

A large class of types of poset congruence are based upon that of a lattice congruence. The idea is that the notion of a lattice congruence is clearly the right one for lattices, and so poset congruences should be defined by extending this notion to posets in a natural way. However, such an extension is clearly not unique, so there exist multiple different notions which operate in this way, or in similar ways.


\subsection{III-congruences}\label{sect:not:iii}

Shum, Zhu, and Kehayopulu introduce the notion of III-congruences and III-homomorph\-isms on posets in order to find a notion of somewhere in between the weak order congruences of Section~\ref{sect:not:njw} and the order congruences of Section~\ref{sect:not:ord} \cite{szk}. One way of looking at III-congruences is that they preserve meet-semilattice structure rather than lattice structure, and consequently are a weaker notion than the later notions in this section.

\begin{definition}[{\cite[Definition~2.5]{szk}}]\label{def:iii}
An equivalence relation $\equ$ on a poset $P$ is called a \emph{III-congruence} if it is a weak order congruence and
\begin{enumerate}
\item given $p, q, r \in P$ such that $p \meet r$ exists, then $p \leqslant q$ and $q \equ r$ implies that $p \equ (p \meet r)$.\label{op:iii:meet}
\end{enumerate}
\end{definition}

The following fact can then be seen from Proposion~\ref{prop:njw:props}\ref{op:njw:works}.

\begin{proposition}
If $\equ$ is a III-congruence on a poset $P$, then the quotient relation $\qrel$ on $\quot$ is a partial order.
\end{proposition}

\begin{definition}[{\cite[Definition~2.1]{szk}}]
A map $f\colon P \to Q$ is called a \emph{III-homomorphism} if it is strong and satisfies the condition that
\begin{itemize}
\item for all $p, q, p', q' \in P$ such that the meets \[p \meet q, p' \meet q', p \meet q \meet p', p' \meet q' \meet q, p \meet q \meet p' \meet q'\] exist, then we have $f(p \meet q) = f(p' \meet q')$ whenever we have both $f(p) = f(q)$ and $f(p') = f(q')$.
\end{itemize}
\end{definition}

In other words, III-homomorphisms are strong maps which preserve certain meets.

\begin{proposition}[{\cite[Theorem~2.9 and Theorem~2.10]{szk}}]
An equivalence relation $\equ$ on a poset $P$ is a III-congruence if and only if it is the kernel of a III-homomorphism $f \colon P \to Q$.
\end{proposition}

The relation between this framework and meet-semilattices is given by the following result.

\begin{proposition}[{\cite[Theorem~3.3]{szk}}]
Let $f \colon M \to M'$ be a surjective map between two meet-semilattices. Then $f$ is a meet-semilattice homomorphism if and only if $f$ is a III-homomorphism.
\end{proposition}

\subsection{\texorpdfstring{$w$}{w}-stable congruences}\label{sect:not:wsc}

Halas introduces the notion of a $w$-stable congruence motivated by the following desiderata \cite{halas}.
\begin{itemize}
\item There should be a well-defined notion of subobject.
\item The intersection of two congruences should be a congruence.
\item The notion should coincide with lattice congruences for lattices.
\end{itemize}

$w$-stable congruences are \textit{only} defined as the kernels of $w$-stable morphisms, rather than also admitting a direct definition in terms of necessary and sufficient conditions on the equivalence relation. Given a poset $P$, we denote by $\maco{P}$ the sublattice of the Dedekind--MacNeille completion $\mac{P}$ generated by the set $\spt{\ls{P}{p} \st p \in P}$. Whilst this is a lattice, it is not always a complete lattice, and so is generally a proper subposet of the Dedekind--MacNeille completion. There is a natural map $\macom{P} \colon P \to \maco{P}$ given by $\macom{P}(p) = \ls{P}{p}$.

\begin{definition}[{\cite{halas}}]
A map $f \colon P \to Q$ is called \emph{$w$-stable} if there is a lattice homomorphism $f^{\ast}\colon \maco{P} \to \maco{Q}$ such that the diagram \[
\begin{tikzcd}
P \ar[r,"f"] \ar[d,"\macom{P}"] & Q \ar[d,"\macom{Q}"] \\
\maco{P} \ar[r,"f^{\ast}"] & \maco{Q}
\end{tikzcd}
\] commutes. That is, $f^{\ast}\iota_{P} = \iota_{Q}f$. A \emph{$w$-stable congruence} is the kernel of a $w$-stable map.
\end{definition}

Note that the definition implies that a $w$-stable map must be order-preserving. It is not generally very easy to check whether a given equivalence relation is $w$-stable.

\begin{proposition}[{\cite[Lemma~9]{halas}}]
If $\equ$ is $w$-stable congruence on a poset $P$, then the transitive closure $\trans{\qrel}$ of the quotient relation $\qrel$ on $\quot$ is a well-defined partial order.
\end{proposition}

\begin{proposition}[\cite{halas}]
Let $L$ be a lattice with $\equ$ an equivalence relation on $L$. Then $\equ$ is a lattice congruence if and only if it is a $w$-stable congruence.
\end{proposition}

$w$-stable congruences on $P$ are precisely the restrictions of lattice congruences on $\maco{P}$. In the following proposition we denote by $\equ^{\ast}|_{P}$ the relation \[\spt{(p, q) \in \equ^{\ast} \st p, q \in P}.\]

\begin{proposition}[{\cite[Lemma~6, Lemma~7]{halas}}]\label{prop:wsc_dmc}
An equivalence relation $\equ$ on $P$ is a $w$-stable congruence if and only if there is a lattice congruence $\equ^{\ast}$ on $\maco{P}$ such that $\equ^{\ast}|_{P} = \equ$.
\end{proposition}

However, the lattice congruence $\equ^{\ast}$ on $\maco{P}$ that restricts to a given $w$-stable congruence $\equ$ on $P$ is not necessarily unique. The set of all $w$-stable congruences on a poset $P$ has good structural properties.

\begin{proposition}[{\cite[Lemma~11]{halas}}]\label{prop:wsc:lattice}
The set of $w$-stable congruences on $P$ forms a complete lattice under inclusion.
\end{proposition}

\begin{remark}
The notion of subobject of a poset mentioned in the desiderata is given by subsets such that the natural inclusion is a $w$-stable map. Such subsets can also be specified by order-theoretic conditions \cite{halas}.
\end{remark}

\subsection{Order congruences}\label{sect:not:ord}

A stronger notion of poset congruence which also derives from lattice congruences is as follows. These congruences possess a good deal of structure, certainly in the finite case (Proposition~\ref{prop:reading}).

\begin{definition}[{\cite[Definition 2]{cs_cong}}]\label{def:cs}
An equivalence $\equ$ on a poset $P$ is called an \emph{order congruence} if
\begin{enumerate}
\item $[p]$ is a convex subset of $P$ for each $p \in P$;\label{op:cs:convex}
\item for each $q, r \in [p]$ there exist $s,t \in [p]$ such that $s \leqslant q \leqslant t$ and $s \leqslant r \leqslant t$;\label{op:cs:meet_join}
\item if $u \leqslant p$, $u \leqslant q$ and $u \equ p$ there exists $t \in P$ with $p \leqslant t$, $q \leqslant t$ and $q \equ t$; 
\item dually, if $p \leqslant v$, $q \leqslant v$ and $v \equ q$ then there exists $s \in P$ with $s \leqslant p$, $s \leqslant q$ and $p \equ s$.\label{op:cs:long}
\end{enumerate}
\end{definition}

\begin{remark}
The definition of a congruence in \cite[Definition~2]{cs_cong} includes the caveat that $P \times P$ is a congruence for any poset $P$, but we exclude this. It is natural to believe that one ought to be able to quotient a poset $P$ by the equivalence relation $P \times P$ to obtain the poset with one element. However, one can accommodate this is with a more permissive notion of congruence altogether, such as those from Section~\ref{sect:not:gen}.
\end{remark}

In the case where $P$ is a finite poset, order congruences admit the following neat description due to Reading \cite{reading_order}.

\begin{proposition}[{\cite[Section 5]{reading_order}}]\label{prop:reading}
Let $P$ be a finite poset with an equivalence relation $\equ$ defined on the elements of $P$. The equivalence relation $\equ$ is an \emph{order congruence} if and only if:
\begin{enumerate}
\item Every equivalence class is an interval.\label{op:read:int}
\item The projection $\pi_{\downarrow} \colon P \to P$, mapping each element $p$ of $P$ to the minimal element in $[p]$, is order preserving.\label{op:read:down}
\item The projection $\pi^{\uparrow} \colon P \to P$, mapping each element $p$ of $P$ to the maximal element in $[p]$, is order-preserving.\label{op:read:up}
\end{enumerate}
\end{proposition}

This notion of congruence gives well-defined quotient posets.

\begin{proposition}
If $\equ$ is an order congruence on a poset $P$, then $\quot$ is a poset.
\end{proposition}

The morphism perspective on order congruences is given as follows. The corresponding morphisms are those that preserve upper and lower bounds.

\begin{definition}[{\cite{cs_cong,reading_order}}]\label{def:order_morph}
A map $f\colon P \to Q$ is an \emph{order morphism} if, for all $p, p' \in P$, we have that \[f(\ls{P}{p, p'}) = \ls{f(P)}{f(p), f(p')}\] and \[f(\us{P}{p, p')} = \us{f(P)}{f(p), f(p')}.\]
\end{definition}

\begin{proposition}[{\cite[Theorem~3]{cs_cong}}]\label{prop:ord:morph_pers}
Let $\equ$ be an equivalence relation on a poset $P$. We have that $\equ$ is an order congruence if and only if it is the kernel of an order morphism $f \colon P \to Q$.
\end{proposition}

\begin{remark}\label{rmk:cs_cong}
Chajda and Sn\'{a}\v{s}el define LU morphisms, whose kernels correspond to their version of order congruences \cite{cs_cong}. A map is automatically an LU morphism if its image has size one, corresponding to the case where the congruence on the poset $P$ is given by $P \times P$. LU morphisms are also required to be surjective, but this does not change the class of congruences obtained.
\end{remark}

\begin{remark}
In fact, upper bounds and lower bounds of finite sets of arbitrary size are preserved by order morphisms, rather than only sets of size two, as in Definition~\ref{def:order_morph}. This is shown in \cite[Lemma~1]{halas}. One can envisage a `completed' version of an order congruence where upper and lower bounds are preserved for sets of arbitrary size, not only finite ones. Indeed, it is shown in \cite[Proposition~2.3]{dirrt} that an analogous description to Proposition~\ref{prop:reading} holds for complete lattice congruences on a complete lattice.
\end{remark}

Order congruences extend the notion of a lattice congruence.

\begin{proposition}
Let $L$ be a lattice with $\equ$ an equivalence relation on $L$. Then $\equ$ is a lattice congruence if and only if it is an order congruence.
\end{proposition}

Order congruences also have a nice interpretation in terms of the Dedekind--MacNeille completion of $P$. Given a finite lattice $L$ with a subposet $P$, a lattice congruence $\equ$ on $L$ \emph{restricts exactly to} $P$ if for every congruence class $[p, q]$ of $\equ$, we either have $p, q \in P$ or $[p, q] \cap P = \emptyset$ \cite{reading_order}.

\begin{proposition}[{\cite[Theorem~8]{reading_order}}]
Let $P$ be a finite poset with Dedekind--MacNeille completion $\mac{P}$, and let $\equ$ be an equivalence relation on $P$. Then $\equ$ is an order congruence on $P$ if and only if there is a lattice congruence $\mac{\equ}$ on $\mac{P}$ which restricts exactly to $P$ such that $\mac{\equ}|_{P} = \equ$, in which case we have that
\begin{enumerate}[label=\textup{(}\arabic*\textup{)}]
\item $\mac{\equ}$ is the unique congruence on $\mac{P}$ which restricts exactly to $\equ$, and
\item the completion $\mac{P/\equ}$ is naturally isomorphic to $\mac{P}/\mac{\equ}$.
\end{enumerate}
\end{proposition}

Compare this to the analogous result for $w$-stable congruences in Proposition~\ref{prop:wsc_dmc}. This gives an idea of the difference between the two types of congruence. Another important difference is that $w$-stable congruences are closed under intersections by Proposition~\ref{prop:wsc:lattice}, while order congruences are not \cite[p.197]{halas}.

\begin{remark}\label{rmk:order_subposet}
Given a finite poset $P$ and an order congruence $\equ$ on $P$, then it can be seen that $\quot$ can be realised as a subposet of $P$, either by sending each equivalence class to its minimal element, or by sending each equivalence class to its maximal element. That this indeed is an embedding of $\quot$ in $P$ follows from Proposition~\ref{prop:reading}. This is one of the features that makes order congruences particularly nice to work with in the finite case.
\end{remark}

\begin{figure}
\caption{A quotient by a lattice congruence which is a subposet but not a sublattice}\label{fig:subpos_not_sublat}
\[
\begin{tikzpicture}[
	frames/.style args = {#1/#2}{minimum height=#1,
               minimum width=#2+\pgfkeysvalueof{/pgf/minimum height},
               draw, rounded corners=3mm, fill=yellow, opacity=0.6,
               sloped},
               ]

	\begin{pgfonlayer}{foreground}
\node(000) at (0,0) {$x_{00}$};
\node(010) at (-1,1) {$x_{10}$};
\node(001) at (1,1) {$x_{01}$};
\node(011) at (0,2) {$x_{11}$};
\node(100) at (0,3) {$y_{00}$};
\node(110) at (-1,4) {$y_{10}$};
\node(101) at (1,4) {$y_{01}$};
\node(111) at (0,5) {$y_{11}$};

\draw[->] (000) -- (001);
\draw[->] (000) -- (010);
\draw[->] (001) -- (011);
\draw[->] (010) -- (011);
\draw[->] (011) -- (100);
\draw[->] (100) -- (101);
\draw[->] (100) -- (110);
\draw[->] (101) -- (111);
\draw[->] (110) -- (111);
	\end{pgfonlayer}

    \begin{pgfonlayer}{main}
\path   let \p1 = ($(011.center)-(100.center)$),
            \n1 = {veclen(\y1,\x1)} in
            (011) -- 
            node[frames=7mm/\n1] {} (100);
    \end{pgfonlayer}

\end{tikzpicture}
\]
\end{figure}
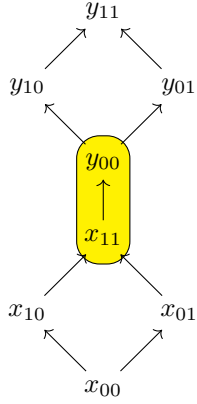

\begin{example}\label{ex:subpos_not_sublat}
It follows from Remark~\ref{rmk:order_subposet} that if $L$ is a finite lattice with a lattice congruence $\equ$ then $L/\equ$ is a subposet of $L$. However, it may not be a sublattice, as remarked in \cite{reading_cambrian}. An example of this is shown in Figure~\ref{fig:subpos_not_sublat}. We illustrate the lattice $L$ by its \emph{Hasse diagram}, which is the graph whose vertices are elements of $L$ with arrows $x \to y$ for covering relations $x \lessdot y$. Here $\equ$ has only one equivalence class which is not a singleton, namely $\{x_{11}, y_{00}\}$. Note first that this poset $L$ is a lattice and that the equivalence relation $\equ$ shown is a lattice congruence. There are two possible ways of embedding $\iota\colon L/\equ \to L$, depending upon whether $\iota([x_{11}]) = x_{11}$ or $\iota([x_{11}]) = y_{00}$. Neither of these embeddings realises $L/\equ$ as a sublattice of $L$. If $\iota([x_{11}]) = x_{11}$, then $y_{10} \meet y_{01} = y_{00} \neq x_{11}$. On the other hand, if $\iota([x_{11}]) = y_{00}$, then $x_{10} \join x_{01} = x_{11} \neq y_{00}$.
\end{example}

\subsection{Haviar--Lihov\'a congruences}\label{sect:not:hl}

Haviar and Lihov\'a \cite{hl_vop} introduce concepts of congruences and homomorphisms of posets which try not only to be consistent with the corresponding notions for lattices, but also those for multilattices, which were introduced in \cite{benado}. A poset $P$ is a \emph{multilattice} if for all $p, q \in P$, and $u \in \us{P}{p, q}$, there exists an element $\overline{u} \in \us{P}{p, q}$ such that $\overline{u} \leqslant u$ and for any $u' \in \us{P}{p, q}$ with $u' \leqslant \overline{u}$, we have $u' = \overline{u}$, with the dual condition holding for $\ls{P}{p, q}$. In other words, a poset is a multilattice if every upper bound of a pair of elements is greater than a minimal upper bound for the pair of elements, with the dual condition holding for lower bounds. Note that finite posets are automatically multilattices.

In order to introduce the notion of a Haviar--Lihov\'a congruence, we need the following notions.

\begin{definition}[{\cite[Definition~2.1]{hl_vop}}]
Let $P$ be a poset with $p, q \in P$. A subset $S \subset \us{P}{p, q}$ is called a \emph{supremum set} or \emph{sup-set} of $p$ and $q$ if the following conditions hold.
\begin{enumerate}
\item For each $u \in \us{P}{p, q}$, there exists $s \in S$ with $s \leqslant u$.
\item If $u \in \us{P}{p, q}$, $s \in S$ with $u \leqslant s$, then $u \in S$.
\end{enumerate}
\emph{Infimum sets}, or \emph{inf-sets} are defined dually.
\end{definition}

In other words, a sup-set is a set of minimal upper bounds for $a$ and~$b$.

\begin{definition}[{\cite[Lemma~4.2]{hl_vop}}]
Let $P$ be a poset. An equivalence relation $\equ$ on $P$ is a \emph{Haviar--Lihov\'a congruence} if and only if it satisfies the following conditions.
\begin{enumerate}
\item All $\equ$-equivalence classes are convex subsets of $P$.
\item If $p, p', q \in P$ are such that $p' \equ p$ and $p \leqslant q$, then $\us{P}{p', q} \neq \emptyset$ and there exists a sup-set $S'$ of $p'$ and $q$ such that $S' \subseteq [q]$.\label{op:hl:sup}
\item Dually, if $p, q, q' \in P$ are such that $p \leqslant q$ and $q \equ q'$, then $\ls{P}{p, q} \neq \emptyset$ and there exists an inf-set $I'$ of $p$ and $q'$ such that $I' \subseteq [p]$.
\end{enumerate}
\end{definition}

The intuition for this definition, in terms of \eqref{op:hl:sup}, is that if $p' \equ p$ and $p \leqslant q$, then the sup-set of $[p']$ and $[q]$ in $\quot$ is $[q]$, so there should be a sup-set $S'$ of $p'$ and $q$ in $P$ which is contained in $[q]$.

\begin{theorem}[{\cite[Theorem~4.8]{hl_vop}}]
Let $P$ be a poset. The intersection of finitely many Haviar--Lihov\'a congruences on $P$ is also a Haviar--Lihov\'a congruence.
\end{theorem}

The morphism perspective on Haviar--Lihov\'a congruences is given as follows.

\begin{definition}[{\cite[Definition~3.1]{hl_vop}}]
Let $P$ and $Q$ be posets with $f \colon P \to Q$ a map. This is called a \emph{Haviar--Lihov\'a homomorphism} if the following conditions are satisfied for all $a, b \in P$.
\begin{enumerate}
\item For each sup-set $S$ of $p, q$, there exists a sup-set $T$ of $f(p)$ and $f(q)$ such that $T \subseteq f(S)$.
\item For each sup-set $T$ of $f(p)$ and $f(q)$, there exists a sup-set $S$ of $p, q$ with $f(S) \subseteq T$.
\item The dual conditions for inf-sets also hold.
\end{enumerate}
\end{definition}

Haviar--Lihov\'a homomorphisms are automatically order-preserving by \cite[Lemma~3.3]{hl_vop}. 

\begin{proposition}[{\cite[Theorem~4.7]{hl_vop}}]\label{prop:hl:morph_pers}
Let $\equ$ be an equivalence relation on a poset $P$. We have that $\equ$ is a Haviar--Lihov\'a congruence if and only if it is the kernel of a Haviar--Lihov\'a morphism $f \colon P \to Q$.
\end{proposition}

The following theorem is the main motivation for Haviar--Lihov\'a congruences.

\begin{theorem}[{\cite[Theorem~3.8]{hl_vop}}]
If $f \colon P \to Q$ is a surjective Haviar--Lihov\'a homomorphism, then
\begin{enumerate}[label=\textup{(}\arabic*\textup{)}]
\item if $P$ is a multilattice, then $Q$ is a multilattice; and
\item if $P$ is a lattice, then $Q$ is a lattice.
\end{enumerate}
\end{theorem}

\begin{proposition}[{\cite[Lemma~3.7,Theorem~4.7]{hl_vop}}]
Let $L$ be a lattice with $\equ$ an equivalence relation on $L$. Then $\equ$ is a lattice congruence if and only if it is a Haviar--Lihov\'a congruence.
\end{proposition}

\begin{remark}
Haviar and Lihov\'a also use their notion of homomorphism to define substructures and varieties of posets. A variety of posets is a class of posets closed under particular operations, such as taking certain homomorphic images and subposets, and taking direct products. This is inspired by Birkhoff's work on varieties of universal algebras \cite{birkhoff_var}. Here a variety of universal algebras is a class of universal algebras possessing certain operations which satisfy certain equations. Birkhoff's Theorem says that a class is a variety of universal algebras if and only if it is closed under homomorphic images, subalgebras, and arbitrary products.
\end{remark}

\section{Further notions}\label{sect:not:misc}

In this section we consider types of poset congruence which do not fit into the groups from Sections~\ref{sect:not:gen} and~\ref{sect:not:lat}, but which instead have different motivation. These notions of congruence arise naturally in examples.

\subsection{GK congruences}\label{sect:not:gk}

The first three notions of congruence in this section all use the same property, which we call being `upper regular'. It is also natural to consider the dual of this property, which we refer to as being `lower regular'.

\begin{definition}[{\cite[pp.48-9]{bj_res}, \cite{blyth}}]
Let $P$ be a poset with $\equ$ an equivalence relation on $P$.
\begin{enumerate}
\item We say that $\equ$ is \emph{upper regular} if, given $p \leqslant q$ and $p \equ p'$, then there exists $q' \in [q]$ such that $p' \leqslant q'$.
\item We say that $\equ$ is \emph{lower regular} if, given $p \leqslant q$ and $q \equ q'$, then there exists $p' \in [p]$ such that $p' \leqslant q'$.
\end{enumerate}
\end{definition}

It is natural to consider these conditions, since they mean that it makes less difference which equivalence-class representatives one chooses when considering the quotient relation.

\begin{remark}
Our terminology here differs from the original terminology from \cite[pp.48-9]{bj_res}, \cite{blyth}. What we call being `upper regular', they call having the `link property'. What they call being `strongly upper regular', we call being `upper regular' and `compatible'.
\end{remark}

\begin{remark}
The condition of being upper regular is the same as the condition that the Bourbaki group call being `weakly compatible' in the context of quotients of pre-orders \cite[Exercise~2, \S 1, Chapter~III]{bourbaki_tos}. The discussion of quotients of posets and pre-orders in \textit{op.\ cit.} is fairly brief.
\end{remark}

Ganesamoorthy and Karpagaval introduce the following natural notion of congruence.

\begin{definition}[{\cite{gk_cong}}]
Let $P$ be a poset with $\equ$ an equivalence relation on $P$. Then we say that $\equ$ is a \emph{GK congruence} if
\begin{enumerate}
\item $\equ$ is upper regular;\label{op:gk:upper}
\item $\equ$ is lower regular; and\label{op:gk:lower}
\item the $\theta$-equivalence classes are convex.\label{op:gk:convex}
\end{enumerate}
\end{definition}

\begin{proposition}
If $\equ$ is a GK congruence on a poset $P$, then $(\quot, \qrel)$ is a well-defined poset.
\end{proposition}

The authors also give a stronger notion for doubly directed posets, in which equivalent pairs of elements have equivalent upper bounds and equivalent lower bounds.

\subsection{Closure congruences}\label{sect:not:clos}

Closure operators on posets seem first to have been considered by Ore \cite{ore_comb,ore_stud,ore_map}, in the context of the poset of subsets of a given set. This originally had nothing to do with taking quotients of posets, but was rather an abstraction of the operation of taking the closure of a set in a topological space. It seems to have been first recognised by Blyth and Janowitz that the kernels of such morphisms gave congruences on posets.

\begin{definition}[{\cite[Theorem~6.9]{bj_res}, \cite{blyth}}]
An equivalence relation $\equ$ on a poset $P$ is a \emph{closure congruence} if and only if
\begin{enumerate}
\item every $\equ$-equivalence class has a unique maximal element, and
\item $\equ$ is upper regular.
\end{enumerate}
\end{definition}

The morphism perspective on closure congruences is given by closure operators. These are now poset endomorphisms, rather than morphisms between two distinct posets, but we still have that closure congruences are precisely the kernels of closure operators.

\begin{definition}[{\cite[p.9]{bj_res}, \cite{blyth}}]
Given a poset $P$, a \emph{closure operator} is an order-preserving map $f \colon P \to P$ such that for all $p \in P$, we have \[f(f(p)) = f(p) \geqslant p.\]
\end{definition}

\begin{proposition}
Let $\equ$ be an equivalence relation on a poset $P$. We have that $\equ$ is a closure congruence if and only if it is the kernel of a closure operator $f \colon P \to P$.
\end{proposition}

Closure operators are also known as `closure relations', and `closure mappings'. The motivating example for a closure operator is of course the poset of subsets of a topological space, with the map $f$ taking a subset to its closure. One may also, of course, study the duals of closure operators and closure congruences.

For a general closure operator $f$ on a poset $P$, we call an element $p \in P$ \emph{closed} if $f(p) = p$. We then have the following result.

\begin{proposition}
Let $P$ be a poset with $f \colon P \to P$ a closure operator and $\equ = \ker f$ the associated closure congruence. Then $(\quot, \qrel)$ is isomorphic to the subposet of closed elements of $P$.
\end{proposition}

Closure congruences can be characterised through their sets of closed elements.

\begin{proposition}[{\cite{mr_clos}}]
Let $P$ be a poset. Then, a subset $S \subseteq P$ is the subset of closed elements of a closure operator if and only if for any $p \in P$, $\us{P}{p} \cap S$ is non-empty and has a unique minimal element.
\end{proposition}

\begin{theorem}[{\cite[Theorem~1, Corollary~6]{hr_clos}}]
The poset of closure congruences on a finite poset $P$ is a join-sublattice of the lattice of equivalence relations on $P$.
\end{theorem}

Further properties of the lattice of closure congruences were proven in \cite{hr_clos}. Indeed, the literature on closure operators is extensive and a full treatment is beyond the scope of this survey. A useful overview is given in \cite{ronse}.

\subsection{Order-autonomous congruences}\label{sect:not:oac}

The notion of an order-autonomous congruence is motivated by that of a lexicographic sum, which is the inverse construction of the quotient by an order-autonomous congruence.
These notions have not historically been viewed in terms of quotients and congruences, but they nevertheless fit neatly into the framework.

\begin{definition}[{\cite{kelly}}]
A non-empty subset $A$ of a poset $P$ is called \emph{order-autonomous} if for all $p \in P \setminus A$, we have that
\begin{enumerate}
\item if $p \leqslant a$ for some $a \in A$, then $p \leqslant a$ for all $a \in A$; and
\item if $p \geqslant a$ for some $a \in A$, then $p \geqslant a$ for all $a \in A$.
\end{enumerate}

An \emph{order-autonomous congruence} on $P$ is one given by a partition of $P$ into order-autonomous subsets.
\end{definition}

\begin{proposition}
If $\equ$ is an order-autonomous congruence on a poset $P$, then $\quot$ is a poset.
\end{proposition}

Order-autonomous congruences arise naturally from the lexicographic-sum construction, which was introduced in \cite{hh} to investigate fixed-point properties of posets.

\begin{definition}[{\cite[Section~3]{hh}}]
Let $\{P_{q}\}_{q \in Q}$ be a family of posets, where $Q$ is itself a poset.
The \emph{lexicographic sum} of the family of posets is the poset \[\lex{\{P_{q}\}_{q \in Q}} := \spt{(q, p) \st q \in Q, p \in P_{q}}\] with partial order given by \[(q, p) \leqslant (q', p') \iff q < q', \text{ or } q = q' \text{ and } p \leqslant p'.\]
\end{definition}

Indeed, we have the following result.

\begin{theorem}
Given a family of posets $\{P_{q}\}_{q \in Q}$, we have that $\{\{q\} \times P_{q}\}_{q \in Q}$ is a partition of $\lex{\{P_{q}\}_{q \in Q}}$ into order-autonomous subsets.
If $\theta$ is the corresponding order-autonomous congruence, then we have that $\quot \cong Q$.

Conversely, suppose that $\{P_{q}\}_{q \in Q}$ is a partition of a poset $P$ into order-auto\-no\-mous subsets giving a congruence $\equ$.
Endowing the set $Q$ with the structure of a poset using the quotient $Q \cong \quot$, one obtains that $P \cong \lex{\{P_{q}\}_{q \in Q}}$.
\end{theorem}

\subsection{Orbits of automorphism groups}\label{sect:not:orb}

Stanley studies quotients of posets by groups of automorphisms. He in particular studies quotients of so-called Peck posets and shows that these quotients retain some nice properties \cite{stanley_quotient_peck}. Other papers studying examples of quotients of posets by groups of automorphisms include \cite{srinivasan,jordan_necklace}.

\begin{proposition}
If $P$ is a finite poset with $\equ$ the equivalence relation given by the orbits of a group $G$ of automorphisms of $P$, then the quotient relation $\qrel$ on $\quot$ is a poset. In this case, we write $P/G := \quot$.
\end{proposition}

Note that if $P$ is an infinite poset with $\equ$ the equivalence relation given by the orbits of a group of automorphisms, then the quotient relation $\qrel$ on $\quot$ is not necessarily a partial order.

Quotients by groups of automorphisms preserve structure which is not always preserved by quotients of posets. Here a poset is \emph{connected} if its Hasse diagram is a connected graph.

\begin{proposition}
Let $P$ be a finite graded connected poset with rank function $\rho$ and $G$ a group of automorphisms of $P$, with $\equ$ the equivalence relation on $P$ given by the orbits of~$G$.
\begin{enumerate}[label=\textup{(}\arabic*\textup{)}]
\item We have that $\rho(p) = \rho(q)$ whenever $p \equ q$ in $P$.
\item Moreover, $P/G$ is a graded poset.
\end{enumerate}
\end{proposition}

\begin{remark}
Quotients of posets by group actions were studied from a different perspective in \cite{bk_gaop}, conceiving posets as a type of loop-free category. The disadvantage of this approach is that it it sometimes produces quotients which are loop-free categories but not posets, and so we omit it.
\end{remark}

\begin{remark}
It is worth briefly remarking on how the quotient of a poset by a group of automorphisms can be conceived using category theory, as is done in \cite{bk_gaop}. Here one conceives a group $G$ as a category $\{\ast\}_{G}$ with one object and all morphisms isomorphisms. An action of a group $G$ on a poset $P$ is then a functor $F \colon \{\ast\}_{G} \to \posetcat$ such that $F(\ast) = P$. The quotient $P/G$ is then the colimit of the functor $F$. The intuition here is that if $f \colon P \to P/G$ is the canonical map and $g \in G$, so that $F(g) \colon P \to P$, then we must have $f \circ F(g) = f$.

Note that this perspective also applies to infinite posets. In cases where the quotient relation is not a partial order, the quotient is produced by finding the smallest compatible congruence containing the equivalence relation given by the orbits, in the manner of Section~\ref{sect:not:comp}.
\end{remark}

\subsubsection{Quotients of Peck posets by groups of automorphisms}

Stanley considers quotients of Peck posets by groups of automorphisms in \cite{stanley_quotient_peck,stanley-appl-alg-comb} and proves a theorem on properties preserved by these sorts of quotients. Let $P$ be a finite graded poset of rank $n$ with ranks $P_{0}, P_{1}, \dots, P_{n}$, as per the notation in Section~\ref{sect:back:def:graded}. If we let $p_{i} = |P_{i}|$, then $P$ is called \emph{rank-symmetric} if $p_{i} = p_{n - i}$ for all $i$ and \emph{rank-unimodal} if there is some $j$ such that $p_{1} \leqslant p_{2} \leqslant \dots \leqslant p_{j} \geqslant p_{j + 1} \geqslant \dots \geqslant p_{n}$. The poset $P$ is \emph{Sperner} if there is no antichain with more elements than the largest of the $p_{i}$. More generally, the poset $P$ is \emph{$k$-Sperner} if the union of $k$ distinct antichains cannot have more elements than the sum of the $k$ largest $p_{i}$. We have that $P$ is \emph{strongly Sperner} if it is $k$-Sperner for $1 \leqslant k \leqslant n + 1$. The poset $P$ is a \emph{Peck poset} if it is rank-symmetric, rank-unimodal, and strongly Sperner.

One can characterise Peck posets using linear algebra. Namely, if $V_{i}$ is the complex vector space with basis $P_{i}$, then we have the following result.

\begin{proposition}[{\cite[Lemma~1.1]{stanley_weyl}}]\label{prop:peck_crit}
The poset $P$ is Peck if and only if there exist linear transformations $\phi_{i}\colon V_{i} \to V_{i + 1}$, for $0 \leqslant i < n$, satisfying the following conditions.
\begin{enumerate}[label=\textup{(}\arabic*\textup{)}]
\item If $p \in P_{i}$, then \[\phi_{i}(p) = \sum_{\substack{q \in P_{i + 1}\\ q > p}} c_{q}q\] for some $c_{q} \in \mathbb{C}$.
\item For all $0 \leqslant i < \frac{1}{2}n$, the linear transformation \[\phi_{n - (i + 1)} \dots \phi_{i + 1}\phi_{i}\colon V_{i} \to V_{n - i}\] is invertible.
\end{enumerate}
\end{proposition}

A Peck poset $P$ is called \emph{unitary} if in the above linear transformations $\phi_{i}$ one can take $c_{q} = 1$ for all $q$.

\begin{theorem}[{\cite[Theorem~1]{stanley_quotient_peck}}]\label{thm:peck}
Let $P$ be a unitary Peck poset with $G$ a group of automorphisms of $P$. Then the quotient poset $P/G$ is Peck.
\end{theorem}

Note that Theorem~\ref{thm:peck} will not hold for other sorts of congruence, since other quotients will not even preserve the property of being graded. Stanley remarks that one cannot do better than Theorem~\ref{thm:peck}: $P/G$ may not be unitary Peck; and $P/G$ may not be Peck if $P$ is Peck but not unitary Peck \cite{stanley_quotient_peck}.

\subsection{Contraction congruences}\label{sect:not:cont}

Contractions seem to have been introduced in \cite{wagner}, although for finite posets their kernels are the connected compatible partitions from \cite{stanley_tpp}, as observed in \cite{thomas_ejc}.
In \cite{cf_decomp}, this notion was independently rediscovered and applied in category theory, as well as related to so-called `admissible' maps of pre-orders \cite{ffm}.

\begin{definition}[{\cite[Section~1]{wagner}}]
A map $f\colon P \to Q$ is a \emph{contraction} if
\begin{enumerate}
\item it is an order-preserving surjection,
\item the fibre $f^{-1}(q)$ is connected for all $q \in Q$, and
\item for any covering relation $q \lessdot q'$ in $Q$, there exists a covering relation $p \lessdot p'$ in $P$ such that $f(p) = q$ and $f(p') = q'$.\label{op:cont:cov}
\end{enumerate}
One can then define a \emph{contraction congruence} to be the kernel of a contraction.
\end{definition}

The intuition for this definition is that the Hasse diagram of $Q$ is the result of contracting convex connected subgraphs of the Hasse diagram of~$Q$. Since \eqref{op:cont:cov} only concerns covering relations, contractions are not in general strong, so we must take the transitive closure of the quotient relation to obtain a well-defined poset.

\begin{proposition}
If $\equ$ is a contraction congruence on a poset $P$, then $(\quot, \trans{\qrel})$ is a well-defined poset.
\end{proposition}

Indeed, for finite $P$, we have the following characterisation of contraction congruences \cite{thomas_ejc}, coming from \cite{stanley_tpp}.

\begin{proposition}
If $P$ is a finite poset with a relation $\equ$, then $\equ$ is a contraction congruence if and only if it is a compatible congruence with connected equivalence classes.
\end{proposition}

\section{Notions for specific types of posets}\label{sect:not:spec}

In this section we survey congruences that have been introduced for specific types of posets.

\subsection{Kolibiar congruences}\label{sect:not:kol}

Kolibiar introduces the following notion of congruence in order to describe direct product decompositions of posets in terms of equivalence relations \cite{kolibiar}. To use this notion of congruence, it is required that the poset be directed. In \cite{kolibiar_ii}, the framework is generalised to connected ordered sets.

\begin{definition}[{\cite[Definition~2.1]{kolibiar}}]\label{def:kolibiar}
Let $P$ be a poset. An equivalence relation $\equ$ on $P$ will be called a \emph{Kolibiar congruence} on $P$ if the following conditions are satisfied.
\begin{enumerate}
\item For each $p \in P$, we have that $[p]$ is a convex subset of $P$.\label{op:kolibiar:1}
\item If $p, q, r \in P$ with $p \leqslant r$, $q \leqslant r$, and $p \equ q$, then there is $s \in P$ such that $p \leqslant s \leqslant r$, $q \leqslant s$ and $p \equ s$.\label{op:kolibiar:2}
\item If $p, q, r, s \in P$, $r \leqslant p \leqslant s$, $r \leqslant q \leqslant s$ and $r \equ p$, then there is $t \in P$ such that $q \leqslant t \leqslant s$, $p \leqslant t$, and $q \equ t$.\label{op:kolibiar:3}
\item The duals of \eqref{op:kolibiar:2} and \eqref{op:kolibiar:3} also hold.
\end{enumerate}
\end{definition}

\begin{remark}
Note that \cite[Definition~2.1]{kolibiar} does not include the dual of~\eqref{op:kolibiar:2}. However, it is clear from \cite[2.2]{kolibiar} that the dual of \eqref{op:kolibiar:2} is intended to hold. Since it does not appear that the dual of \eqref{op:kolibiar:2} follows from the remaining conditions, we include it explicitly.
\end{remark}

It is not generally very easy to check whether a given equivalence relation is a Kolibiar congruence. The defining properties of a Kolibiar congruence are chosen to make the following result true.

\begin{proposition}
An equivalence relation $\equ$ on a lattice $L$ is a lattice congruence if and only if it is a Kolibiar congruence.
\end{proposition}

Hence, the notion could also go in Section~\ref{sect:not:lat}. However, we put it here since it has Theorem~\ref{thm:kol} as a specific piece of motivation, which, along with the other key properties of the notion, is only shown for directed posets in \cite{kolibiar}. As ever, the easiest way to verify the above proposition is by applying the criterion Lemma~\ref{lem:tech}. Indeed, Definition~\ref{def:kolibiar}\eqref{op:kolibiar:2} roughly corresponds to Lemma~\ref{lem:tech}\ref{op:tech:collapse}, whilst Definition~\ref{def:kolibiar}\eqref{op:kolibiar:3} roughly corresponds to Lemma~\ref{lem:tech}\ref{op:tech:subs}. The other key properties of Kolibiar congruences are as follows.

\begin{proposition}[{\cite[2.6, Theorem~2, Theorem~3]{kolibiar}}]
Let $\equ$ be a Kolibiar congruence on a directed poset $P$.
\begin{enumerate}[label=\textup{(}\arabic*\textup{)}]
\item We have that $(\quot, \qrel)$ is a well-defined poset.
\item The set of Kolibiar congruences on $P$ forms a complete distributive lattice which is a sublattice of the lattice of equivalence relations on~$P$.
\end{enumerate}
\end{proposition}

The principal motivation for Kolibiar congruences is the following theorem.

\begin{theorem}[{\cite[Theorem~7]{kolibiar}}]\label{thm:kol}
Let $P$ be a directed poset. Then there is a bijection between direct product decompositions $P = \prod_{i = 1}^{n}P_{i}$ and families $\spt{\equ_{i} \st i \in [n]}$ of Kolibiar congruences satisfying the following conditions.
\begin{enumerate}[label=\textup{(}\arabic*\textup{)}]
\item $\bigcap_{i = 1}^{n} \equ_{i} = \mathrm{id}$.
\item $\bigvee_{i = 1}^{n} \equ_{i} = P \times P$, where $\vee$ denotes the smallest equivalence relation containing a given set of equivalence relations.
\item Given a set $\{p_{1}, \dots, p_{n}\} \subseteq P$, there exists an element $p \in P$ such that $p\equ_{i}p_{i}$ for all $i \in [n]$.
\end{enumerate}
\end{theorem}

In terms of the direct product decomposition, the equivalence relation $\equ_{i}$ should be thought of as identity in the $i$-th coordinate: $p \equ_{i} q$ if and only if the $i$-th coordinates of $p$ and $q$ are the same. One can think of Kolibiar congruences as congruences that arise in this way.

The morphism perspective on Kolibiar congruences is given by the following definition and result.

\begin{definition}
A map $f \colon P \to Q$ between two directed posets is called a \emph{Kolibiar morphism} if
\begin{enumerate}
\item $f$ is order-preserving;\label{op:kol_map:1}
\item if $p, q, r \in P$ with $p \leqslant r$ and $q \leqslant r$ and $f(p) \leqslant f(q)$, then there is $s \in P$ such that $p \leqslant r$, $q \leqslant s \leqslant r$ and $f(q) = f(s)$;\label{op:kol_map:2}
\item the dual of \eqref{op:kol_map:2} also holds.
\end{enumerate}
\end{definition}

\begin{proposition}[{\cite[Theorem~1]{kolibiar}}]
Let $\equ$ be an equivalence relation on a directed poset $P$. We have that $\equ$ is a Kolibiar congruence if and only if it is the kernel of a Kolibiar morphism $f \colon P \to Q$.
\end{proposition}

\subsection{Homogeneous congruences}\label{sect:not:hom}

Homogeneous quotients were introduced by Hallam and Sagan in \cite{hs-char-poly} to study the characteristic polynomials of posets. See also \cite{hallam_applications}. The definition requires that the poset has a unique minimal element.

\begin{definition}[{\cite[Definition 4]{hs-char-poly}}]\label{def:hs}
Let $P$ be a finite poset with a unique minimal element $\hat{0}$, and let $\equ$ be an equivalence relation on $P$. Then we say that $\equ$ is a \emph{homogeneous congruence} if
\begin{enumerate}
\item $\{\hat{0}\}$ is a $\equ$-equivalence class, and\label{op:hs:min}
\item $\equ$ is upper regular.\label{op:hs:important}
\end{enumerate}
\end{definition}

\begin{proposition}[{\cite{hs-char-poly}}]\label{prop:qh:wd}
If $\equ$ is a homogeneous congruence on a finite poset $P$, then we have that $(\quot, \qrel)$ is a well-defined poset.
\end{proposition}

Note that here it is required that the poset be finite.

\begin{remark}
Note that property~\eqref{op:hs:min} is not required in the proof the quotient is well-defined in Proposition~\ref{prop:qh:wd}. The motivation for this property comes from the fact that characteristic polynomials are only defined for posets with a unique minimal element.
\end{remark}

\subsubsection{M\"obius functions and characteristic polynomials}

We explain more detail on the motivation for homogeneous congruences, which stems from characteristic polynomials of posets. Given a poset $P$ with a unique minimal element $\hat{0}$, recall that the (\emph{one-variable}) \emph{M\"obius function} of $P$ is the function $\mu_{P}\colon P \to \mathbb{Z}$ defined recursively~by \[\sum_{p \leqslant q} \mu_{P}(p) = \delta_{\hat{0},q},\] where $\delta_{\hat{0}, q}$ is the Kronecker delta.

We now suppose that $P$ is graded with rank function $\rho \colon P \to \mathbb{Z}_{\geqslant 0}$. We denote the rank of $P$ by $\rho(P)$. The \emph{characteristic polynomial} of $P$ is the generating function for $\mu$, that is, \[\chi_{P}(t) = \sum_{p \in P}\mu_{P}(p)t^{\rho(P) - \rho(p)}.\] Hallam and Sagan investigate the characteristic polynomials of lattices using homogeneous quotients in \cite{hs-char-poly}. If a homogeneous congruence satisfies certain conditions, then quotienting does not affect the characteristic polynomial.

\begin{proposition}[{\cite[Lemma~6, Corollary~8]{hs-char-poly}}]\label{prop:char_poly}
Let $P$ be a graded poset with a unique minimal element $\hat{0}$ with $\equ$ a homogeneous congruence on $P$. Suppose that, for all $[p] \in \quot$ with $\hat{0} \notin [p]$, we have \[\sum_{q \in \ls{P}{[p]}} \mu_{P}(q) = 0,\] and that the rank function $\rho$ is constant on equivalence classes. Then \[\mu_{\quot}([p]) = \sum_{q \in [p]} \mu_{P}(q)\] and, consequently, \[\chi_{\quot}(t) = \chi_{P}(t).\]
\end{proposition}

By applying Proposition~\ref{prop:char_poly}, one can compute the characteristic polynomials of posets by taking certain homogeneous quotients in order to simplify them.

\begin{remark}
There are other approaches one could take to quotienting posets in such a way that preserves the characteristic polynomial. For instance, one could replace the condition that, for all $[p] \in P/\theta$ with $\hat{0} \notin [p]$, we have that
\begin{equation}\label{eq:hom_old}
\sum_{q \in \ls{P}{[p]}}\mu_{P}(q) = 0,
\end{equation}
with the condition that for all $[p] \in P/\theta$ we have
\begin{equation}\label{eq:hom_new}
\sum_{[q] \in \ls{\quot}{[p]}} \sum_{r \in [q]} \mu_{P}(r) = \indic{[q]}{\hat{0}}.
\end{equation}
Here $\indic{X}{x}$ is the indicator function of $X$, which equals $1$ if $x \in X$ and $0$ otherwise. This condition is equivalent to having $\mu_{\quot}([p]) = \sum_{q \in [p]}\mu_{P}(q)$. The constancy of the rank function on equivalence classes then gives the desired result that $\chi_{\quot}(t) = \chi_{P}(t)$. Note that by Proposition~\ref{prop:grade&comp}, if the rank function $\rho$ is constant on $\equ$-equivalence classes, then $\equ$ is compatible and $\trans{\qrel}$ is a well-defined partial order on $\quot$. Hence, by replacing the assumption~\eqref{eq:hom_old} with the assumption~\eqref{eq:hom_new}, while maintaining the assumption that the rank function is constant on equivalence classes, no additional assumption is then needed on the equivalence relation $\equ$. This then gives a strictly larger set of congruences which preserve the characteristic polynomial.
\end{remark}

In \cite{hallam_applications}, Hallam uses homogeneous quotients to give proofs of classical results on M\"obius functions by induction on the size of the poset.

\subsection{Natural DCPO congruences}\label{sect:not:dcpo}

Mahmoudi, Moghbeli, and Pi\'oro study congruences of posets which are directed complete \cite{mmp_nat,mmp_dcpo}. Here, a poset $P$ is \emph{directed complete} if every directed subposet $D \subseteq P$ has a join $\bigvee D$ in $P$. `Directed complete partial order' is abbreviated `DCPO'. DCPOs are fundamental in Domain Theory, a mathematical foundation of the semantics of programming languages introduced by Scott \cite{scott69,scott70}. See \cite{aj_dt} for an overview.

\begin{definition}[{\cite{scott70}}]
Let $P$ and $Q$ be DCPOs. A \emph{DCPO map} $\phi \colon P \to Q$ is a map such that for each directed subposet $D \subseteq P$, we have that the subposet $\phi(D)$ of $Q$ is directed, and $\phi\left( \bigvee D \right) = \bigvee \phi(D)$.
\end{definition}

DCPO maps are also known as `continuous' or `Scott-continuous'. Note that it follows from the definition that DCPO maps are order-preserving, since a pair of comparable elements in $P$ forms a directed subposet. The fundamental theorem of domain theory is that every DCPO map $P \to P$ has a least fixed point; this least fixed point is then the mathematical counterpart of a recursive definition in a program.

\begin{definition}[\cite{mmp_nat}]
An equivalence relation $\equ$ is a \emph{natural DCPO congruence} if the transitive closure $\trans{\qrel}$ of the quotient relation $\qrel$ on $\quot$ is a DCPO, with the canonical map $P \to \quot$ a DCPO map.
\end{definition}

This sort of definition follows a recipe that can be used to define a notion of congruence for any type of poset with additional structure, such as a DCPO. Namely, a congruence should be an equivalence relation $\equ$ such that the canonical map $P \to \quot$ preserves the additional structure. It ought to be possible to define DCPO congruences in terms of necessary and sufficient conditions on the equivalence relation, however. Note also that for finite DCPOs, natural DCPO congruences coincide with compatible congruences.

\begin{remark}
In their work, Mahmoudi, Mohgbeli, and Pi\'oro also consider `DCPO congruences', which are precisely the kernels of DCPO maps. However, given a DCPO congruence $\equ$ on a DCPO $P$, there is no canonically defined DCPO on $\quot$. Indeed, the transitive closure $\trans{\qrel}$ of the quotient relation $\qrel$ on $\quot$ may not be a DCPO, and there may be several ways of extending this partial order to a DCPO. In general, therefore, `DCPO congruences' are not always natural DCPO congruences, and so natural DCPO congruences are not precisely the kernels of DCPO maps.
\end{remark}

\subsection{Posets with additional algebraic structure}\label{sect:not:alg}

We consider briefly a couple of types of posets with additional algebraic structure. Similar to the case of lattices, congruences on such posets are also required to respect the algebraic structure, in the sense from universal algebra. This situation is considered in general in \cite{krishnan}.

\subsubsection{Hilbert algebras}

Hilbert algebras give another instance, analogous to lattices, where a poset has additional algebraic structure, namely a binary operation $\ast$ obeying certain axioms. One can then consider congruences of Hilbert algebras according to the universal algebra approach. However, unlike lattices, congruences of Hilbert algebras do not always produce well-defined posets \cite{clp_rel}. The framework of Hilbert algebras was used to study congruences of relatively pseudocomplemented posets in \cite{clp_rel}, since relatively pseudocomplemented posets form a subclass of the class of Hilbert algebras. Hilbert algebras arise from the axiomatisation of propositional logic given by Hilbert in \cite{hilbert_logic}, and it seems that they were first put into the framework of posets in \cite{rasiowa}.

\subsubsection{Sectionally pseudocomplemented posets}\label{sect:not:spp}

Another type of posets with additional algebraic structure are sectionally pseudocomplemented posets.
The additional algebraic structure determines the order if the poset has a unique maximal element \cite{clp_spp}. For these posets a congruence on the algebra does not necessarily guarantee that the quotient relation is a partial order. However, this is guaranteed when the poset has a unique maximal element, satisfies the ascending chain condition, and is \textit{strongly} sectionally pseudocomplemented \cite[Theorem~3.2 and Theorem~3.5]{clp_spp}. (Recall that the ascending chain condition holds when there are no infinite strictly ascending chains.) Sectionally pseudocomplemented posets play a role in semantics in a similar way to DCPOs: see \cite{cl_imp} and references therein.

\section{Examples}\label{sect:ex}

We now survey examples of quotient posets that have appeared in the literature. We focus on giving nice examples, rather than representing all different types of congruences. Other examples from the literature include \cite{dhp,hossain_thesis,lw_minor}.

\subsection{Graphs}

Our first example shows how the poset of graphs on a given set of vertices can be constructed as a quotient by an automorphism group. Let $n = \binom{m}{2}$ for some positive integer $m$. We consider the Boolean lattice $B_{n}$ of subsets of $[n]$ as the set of labelled simple graphs on $m$ vertices, as we presently explain. \emph{Simple} graph here means that there is at most one edge between any two vertices. Identify each of the elements of $[n]$ with a different unordered pair $\{i, j\}$ of two distinct elements from $[m]$, corresponding to an edge from $i$ to~$j$. The set $B_{n}$ then naturally corresponds to the set of labelled simple graphs on $m$ vertices, with the order relation given by edge inclusion.

Let $\binom{[m]}{2}$ denote the set of subsets of $[m]$ of size $2$, and let $G$ be the symmetric group $S_{m}$. This acts by permuting the $m$ points, which induces an action on the set of edges $\binom{[m]}{2}$, and hence on the poset $B_{n}$. This action on the poset $B_{n}$ simply relabels the vertices of the graph, and so orbits correspond to isomorphic graphs on the $m$ vertices. Hence, the quotient poset $B_{n}/G$ is the subgraph order on the set of non-isomorphic simple graphs on $m$ vertices.  Applying Theorem~\ref{thm:peck} gives the following result.

\begin{theorem}[{\cite{stanley-appl-alg-comb,stanley_quotient_peck}}]
The poset of non-isomorphic simple graphs on $m$ vertices with respect to subgraph inclusion is Peck.
\end{theorem}

\begin{figure}
\caption{The subgraph ordering on simple graphs on four vertices}\label{fig:graph_four}
\[
\begin{tikzpicture}[xscale=1.5,yscale=1.2]

\node (0) at (0,0) {$\bullet$};
\node (1) at (0,1) {$\bullet$};
\node (20) at (-1,2) {$\bullet$};
\node (21) at (1,2) {$\bullet$};
\node (30) at (-2,3) {$\bullet$};
\node (31) at (-1,3) {$\bullet$};
\node (32) at (0,3) {$\bullet$};
\node (40) at (-1,4) {$\bullet$};
\node (41) at (1,4) {$\bullet$};
\node (5) at (0,5) {$\bullet$};
\node (6) at (0,6) {$\bullet$};

\draw[->] (0) -- (1);
\draw[->] (1) -- (20);
\draw[->] (1) -- (21);
\draw[->] (20) -- (30);
\draw[->] (20) -- (31);
\draw[->] (20) -- (32);
\draw[->] (21) -- (32);
\draw[->] (30) -- (40);
\draw[->] (31) -- (40);
\draw[->] (32) -- (40);
\draw[->] (32) -- (41);
\draw[->] (40) -- (5);
\draw[->] (41) -- (5);
\draw[->] (5) -- (6);


\begin{scope}[scale=0.25,shift={(-1.25,-2)},font=\tiny]

\coordinate (1) at (0,0);
\coordinate (2) at (1,0);
\coordinate (3) at (2,0);
\coordinate (4) at (3,0);

\node at (1) {$\bullet$};
\node at (2) {$\bullet$};
\node at (3) {$\bullet$};
\node at (4) {$\bullet$};

\end{scope}


\begin{scope}[scale=0.25,shift={(1,3.5)},font=\tiny]

\coordinate (1) at (0,0);
\coordinate (2) at (1,0);
\coordinate (3) at (2,0);
\coordinate (4) at (3,0);

\draw (3) -- (4);

\node at (1) {$\bullet$};
\node at (2) {$\bullet$};
\node at (3) {$\bullet$};
\node at (4) {$\bullet$};

\end{scope}


\begin{scope}[scale=0.25,shift={(-7,8)},font=\tiny]

\coordinate (1) at (0,0);
\coordinate (2) at (1,0);
\coordinate (3) at (2,0);
\coordinate (4) at (1,-1);

\draw (1) -- (2);
\draw (2) -- (3);

\node at (1) {$\bullet$};
\node at (2) {$\bullet$};
\node at (3) {$\bullet$};
\node at (4) {$\bullet$};

\end{scope}


\begin{scope}[scale=0.25,shift={(5,7.5)},font=\tiny]

\coordinate (1) at (0,0);
\coordinate (2) at (1,0);
\coordinate (3) at (0,1);
\coordinate (4) at (1,1);

\draw (1) -- (2);
\draw (3) -- (4);

\node at (1) {$\bullet$};
\node at (2) {$\bullet$};
\node at (3) {$\bullet$};
\node at (4) {$\bullet$};

\end{scope}


\begin{scope}[scale=0.25,shift={(-10,12)},font=\tiny]

\coordinate (1) at (0,0);
\coordinate (2) at (1,0);
\coordinate (3) at (0.5,1);
\coordinate (4) at (0.5,-1);

\draw (1) -- (2);
\draw (2) -- (3);
\draw (3) -- (1);

\node at (1) {$\bullet$};
\node at (2) {$\bullet$};
\node at (3) {$\bullet$};
\node at (4) {$\bullet$};

\end{scope}


\begin{scope}[scale=0.25,shift={(-2.5,12)},font=\tiny]

\coordinate (1) at (-90:1);
\coordinate (2) at (0,0);
\coordinate (3) at (55:1);
\coordinate (4) at (125:1);

\draw (1) -- (2);
\draw (2) -- (3);
\draw (2) -- (4);

\node at (1) {$\bullet$};
\node at (2) {$\bullet$};
\node at (3) {$\bullet$};
\node at (4) {$\bullet$};

\end{scope}


\begin{scope}[scale=0.25,shift={(1.5,12)},font=\tiny]

\coordinate (1) at (0,0);
\coordinate (2) at (1,0);
\coordinate (3) at (2,0);
\coordinate (4) at (3,0);

\draw (1) -- (2);
\draw (2) -- (3);
\draw (3) -- (4);

\node at (1) {$\bullet$};
\node at (2) {$\bullet$};
\node at (3) {$\bullet$};
\node at (4) {$\bullet$};

\end{scope}


\begin{scope}[scale=0.25,shift={(-6,16)},font=\tiny]

\coordinate (1) at (0,0);
\coordinate (2) at (1,0);
\coordinate (3) at (0.5,1);
\coordinate (4) at (-0.5,1);

\draw (1) -- (2);
\draw (2) -- (3);
\draw (3) -- (1);
\draw (3) -- (4);

\node at (1) {$\bullet$};
\node at (2) {$\bullet$};
\node at (3) {$\bullet$};
\node at (4) {$\bullet$};

\end{scope}


\begin{scope}[scale=0.25,shift={(5,15.5)},font=\tiny]

\coordinate (1) at (0,0);
\coordinate (2) at (1,0);
\coordinate (3) at (0,1);
\coordinate (4) at (1,1);

\draw (1) -- (2);
\draw (3) -- (4);
\draw (2) -- (4);
\draw (1) -- (3);

\node at (1) {$\bullet$};
\node at (2) {$\bullet$};
\node at (3) {$\bullet$};
\node at (4) {$\bullet$};

\end{scope}


\begin{scope}[scale=0.25,shift={(1.5,19.5)},font=\tiny]

\coordinate (1) at (0,0);
\coordinate (2) at (1,0);
\coordinate (3) at (0,1);
\coordinate (4) at (1,1);

\draw (1) -- (2);
\draw (3) -- (4);
\draw (2) -- (4);
\draw (1) -- (3);
\draw (2) -- (3);

\node at (1) {$\bullet$};
\node at (2) {$\bullet$};
\node at (3) {$\bullet$};
\node at (4) {$\bullet$};

\end{scope}


\begin{scope}[scale=0.25,shift={(-2,23.5)},font=\tiny]

\coordinate (1) at (0,0);
\coordinate (2) at (1,0);
\coordinate (3) at (0,1);
\coordinate (4) at (1,1);

\draw (1) -- (2);
\draw (3) -- (4);
\draw (2) -- (4);
\draw (1) -- (3);
\draw (2) -- (3);
\draw (1) -- (4);

\node at (1) {$\bullet$};
\node at (2) {$\bullet$};
\node at (3) {$\bullet$};
\node at (4) {$\bullet$};

\end{scope}

\end{tikzpicture}
\]
\end{figure}
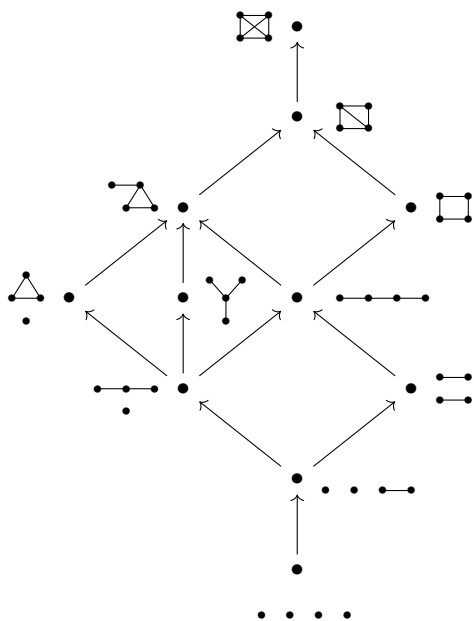

An example of such a poset is shown in Figure~\ref{fig:graph_four}, \emph{cf.} \cite[Figure~4]{stanley-appl-alg-comb}.

\subsection{Young diagrams}\label{sect:ex:young}

We now similarly examine how the poset of Young diagrams contained in a given rectangle may be constructed as a quotient by a group of automorphisms. We consider the Boolean lattice $B_{mn}$ and think of the underlying set of $mn$ elements as a rectangular array of cells with $m$ rows and $n$ columns. The wreath product $G = S_{m} \wr S_{n}$ permutes the $n$ cells within each row independently, and permutes the $m$ rows by interchanging them. The group $G$ thus has $|G| = (n!)^{m}m!$ elements.

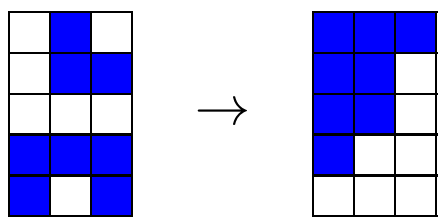
\begin{figure}
\caption{The Young diagram in the orbit of a set of cells}\label{fig:young_orbit}
\[
\begin{tikzpicture}

\node at (0,0) {
\begin{ytableau}
*(white) & *(blue) & *(white) \\
*(white) & *(blue) & *(blue) \\
*(white) & *(white) & *(white) \\
*(blue) & *(blue) & *(blue) \\
*(blue) & *(white) & *(blue)
\end{ytableau}
};

\node at (2,0) {\Huge $\to$};

\node at (4,0) {
\begin{ytableau}
*(blue) & *(blue) & *(blue) \\
*(blue) & *(blue) & *(white) \\
*(blue) & *(blue) & *(white) \\
*(blue) & *(white) & *(white) \\
*(white) & *(white) & *(white)
\end{ytableau}
};

\end{tikzpicture}
\]
\end{figure}

Given a set of cells $T \in B_{mn}$, there is a canonical representative of its orbit under $G$, which is obtained as follows. First, the cells of the rows of the array are permuted to move the cells of $T$ as far left as possible. The rows of the array are then interchanged so that the number of cells in each row decreases as one goes down the array. A finite collection of cells arranged in left-justified rows of decreasing length is called a \emph{Young diagram}. These correspond to integer partitions, with the entries of the partition corresponding to the lengths of the rows. See Figure~\ref{fig:young_orbit} for an example.

\begin{figure}
\caption{The poset $L(2, 3)$}\label{fig:young_23}
\[
\begin{tikzpicture}[xscale=1.5,yscale=1.2]

\node (0) at (0,0) {$\bullet$};
\node (1) at (0,1) {$\bullet$};
\node (20) at (-1,2) {$\bullet$};
\node (21) at (1,2) {$\bullet$};
\node (30) at (-2,3) {$\bullet$};
\node (31) at (0,3) {$\bullet$};
\node (40) at (-1,4) {$\bullet$};
\node (41) at (1,4) {$\bullet$};
\node (5) at (0,5) {$\bullet$};
\node (6) at (0,6) {$\bullet$};

\draw[->] (0) -- (1);
\draw[->] (1) -- (20);
\draw[->] (1) -- (21);
\draw[->] (20) -- (30);
\draw[->] (20) -- (31);
\draw[->] (21) -- (31);
\draw[->] (30) -- (40);
\draw[->] (31) -- (40);
\draw[->] (31) -- (41);
\draw[->] (40) -- (5);
\draw[->] (41) -- (5);
\draw[->] (5) -- (6);

\ytableausetup{smalltableaux}


\begin{scope}[shift={(0,-0.3)}]

\node at (0,0) {$\varnothing$};

\end{scope}


\begin{scope}[shift={(0.4,1)}]

\node at (0,0) {
\begin{ytableau}
*(blue)
\end{ytableau}
};

\end{scope}


\begin{scope}[shift={(-1.5,2)}]

\node at (0,0) {
\begin{ytableau}
*(blue) & *(blue)
\end{ytableau}
};

\end{scope}


\begin{scope}[shift={(1.25,2)}]

\node at (0,0) {
\begin{ytableau}
*(blue) \\
*(blue)
\end{ytableau}
};

\end{scope}


\begin{scope}[shift={(-2.5,3)}]

\node at (0,0) {
\begin{ytableau}
*(blue) & *(blue) & *(blue)
\end{ytableau}
};

\end{scope}


\begin{scope}[shift={(0.75,3)}]

\node at (0,0) {
\begin{ytableau}
*(blue) & *(blue) \\
*(blue)
\end{ytableau}
};

\end{scope}


\begin{scope}[shift={(-1.25,4.25)}]

\node at (0,0) {
\begin{ytableau}
*(blue) & *(blue) & *(blue) \\
*(blue)
\end{ytableau}
};

\end{scope}


\begin{scope}[shift={(1.4,4)}]

\node at (0,0) {
\begin{ytableau}
*(blue) & *(blue) \\
*(blue) & *(blue)
\end{ytableau}
};

\end{scope}


\begin{scope}[shift={(0.6,5.25)}]

\node at (0,0) {
\begin{ytableau}
*(blue) & *(blue) & *(blue) \\
*(blue) & *(blue)
\end{ytableau}
};

\end{scope}


\begin{scope}[shift={(0,6.5)}]

\node at (0,0) {
\begin{ytableau}
*(blue) & *(blue) & *(blue) \\
*(blue) & *(blue) & *(blue)
\end{ytableau}
};

\end{scope}

\end{tikzpicture}
\]
\end{figure}
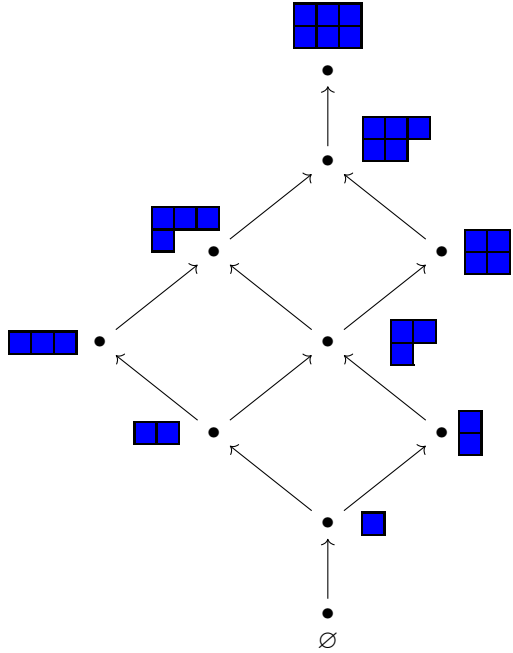

As before, we consider the quotient poset $B_{mn}/G$. Each orbit of the action of $G$ on $B_{mn}$ contains precisely one Young diagram, so this poset is the poset $L(m, n)$ of Young diagrams contained in an $m \times n$ rectangle ordered by inclusion. The poset $L(m, n)$ is in fact a distributive lattice. An example of one of these posets is shown in Figure~\ref{fig:young_23}. By applying Theorem~\ref{thm:peck}, we obtain the following result.

\begin{theorem}[{\cite{stanley-appl-alg-comb,stanley_quotient_peck}}]\label{thm:young}
The poset $L(m, n)$ is Peck.
\end{theorem}

Note that, despite the fact that both $B_{mn}$ and $L(m, n)$ are lattices, the congruence on the former which gives the latter is not a lattice congruence, due to Proposition~\ref{prop:orb_int_ord=triv}. Hence there are interesting quotients from lattices to lattices which are not given by lattice congruences.

\subsection{The poset of conjugacy classes of subgroups}

A well-studied example of a quotient of a poset by a group of automorphisms is given by the poset of conjugacy classes of subgroups of a particular group. Indeed, let $G$ be a finite group with $\Lambda(G)$ its lattice of subgroups. There is a natural action of $G$ on $\Lambda(G)$ by conjugacy. Indeed, given a subgroup $H \leqslant G$ and an element $g \in G$, the action of $g$ on $H$ is given by $H^g := g^{-1}Hg$. The \emph{poset of conjugacy classes of subgroups of $G$} is then the quotient poset $\mathcal{C}(G) := \Lambda(G)/G$. This poset is sometimes called the \emph{frame} of~$G$.

The frame of a group was first studied in \cite{hio_moebius}, where the following theorem on the M\"obius functions of $\Lambda(G)$ and $\mathcal{C}(G)$ was proven. Here $[G,G]$ denotes the subgroup of commutators of $G$.

\begin{theorem}[{\cite[Theorem~7.2]{hio_moebius}}]
If $G$ is solvable, then \[\mu_{\Lambda(G)}(G) = \mu_{\mathcal{C}(G)}(G)|[G,G]|.\]
\end{theorem}

Another connection between the solvability of a group and its frame was shown in \cite{fumagalli}. Recall here that if $P$ is a poset with maximal element $\hat{1}$, then $p \in P$ is a \emph{coatom} if $p$ is covered by $\hat{1}$. Recall also that a lattice $L$ is \emph{modular} if for every $x, y, z \in L$ such that $x \leqslant z$, we have that \[(x \join y) \meet z = x \join (y \meet z).\]

\begin{theorem}[{\cite{fumagalli}}]
A finite group $G$ is solvable if and only if every collection of coatoms of $\mathcal{C}(G)$ has a well-defined meet and the poset consisting of all such  meets is a modular lattice.
\end{theorem}

There is a fair amount of literature on the frame of a group. The homotopy type of the order complex of $\mathcal{C}(G)$ was studied in \cite{welker}. See also \cite{tarnauceanu,mainardis,bcr} for other results on this poset. There are variations on the frame, such as looking at: the poset of conjugacy classes of a restricted set of subgroups \cite{ww}; the quotient of the poset of subgroups by identifying all isomorphic subgroups, rather than just conjugate ones \cite{tarnauceanu_iso}; and the quotient of the subgroup lattice of $G$ by a group other than $G$ itself \cite{mainardis_p}.

\subsection{Permutations to binary trees}\label{sect:ex:perms->trees}

There is a natural order on the set of binary trees with a given number of vertices known as the \emph{Tamari lattice}. This poset can be expressed as a quotient of the weak Bruhat order on the symmetric group, giving a map from permutations to binary trees. This map appears in many guises, including as a map from the permutohedron to the associahedron \cite{tonks}, in the context of Hopf algebras \cite{lr_hopf,lr_order,loday_dialgs}, and in Coxeter theory \cite{bw_coxeter,bw_shell_2}. See also \cite[Section~1.3.13]{stanley_enum} and \cite{thomas-bst}. The map moreover appears as the prototype for the theory of Cambrian lattices \cite{reading_cambrian}, as we explain in the next section. 

Our framework for binary trees is based on \cite[Section~1.5, Chapter~1]{tarjan}. A \emph{tree} is an undirected graph that is connected and acyclic. A \emph{rooted tree} is a tree $T$ with a distinguished vertex $r$, called the \emph{root}. Given vertices $v$ and $w$ such that $v$ is on the path from $r$ to $w$, then we say that $v$ is an \emph{ancestor} of $w$ and $w$ is a \emph{descendant} of $v$. If additionally $v \neq w$, then $v$ is a \emph{proper ancestor} of $w$ and $w$ is a \emph{proper descendant} of $v$. If $v$ is a proper ancestor of $w$ with the two vertices adjacent, then we say that $v$ is the \emph{parent} of $w$ and $w$ is a \emph{child} of $v$. A \emph{binary tree} is a rooted tree in which each vertex $v$ either has no children, or has exactly two children, namely its \emph{left child} $l(v)$ and its \emph{right child} $r(v)$. A vertex with two children is \emph{internal} and a vertex with no children is \emph{external}. A binary tree is said to be \emph{of size $n$} if it has $n$ internal nodes. A tree of size $n$ has $n + 1$ external nodes. Given an internal vertex $v$, its \emph{left subtree} is the subtree rooted at its left child and its \emph{right subtree} is the subtree rooted at its right child.

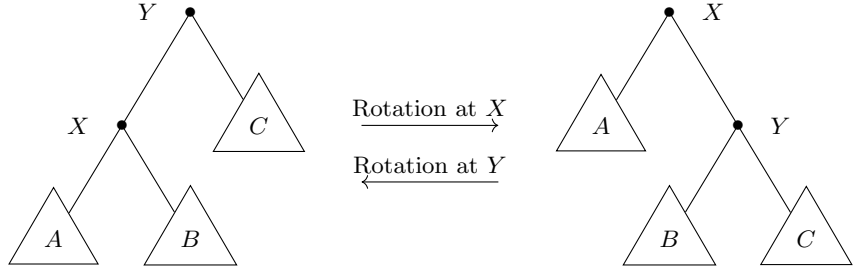
\begin{figure}
\caption{Rotation of binary trees (cf.\ \cite[Figure~1]{stt})}\label{fig:rotation}
\[
\begin{tikzpicture}[scale=0.75,font=\small,xscale=1.2]


\begin{scope}[yscale=2,shift={(-3.5,0)}]

\node (Y) at (0,0) {$\bullet$};
\node at (Y) [left = 1mm of Y]{$Y$};
\node (X) at (-1,-1) {$\bullet$};
\node at (X) [left = 1mm of X]{$X$};
\node [regular polygon,regular polygon sides=3,draw] (A) at (-2,-2) {$A$};
\node [regular polygon,regular polygon sides=3,draw] (B) at (0,-2) {$B$};
\node [regular polygon,regular polygon sides=3,draw] (C) at (1,-1) {$C$};

\draw (A) -- (-1,-1);
\draw (B) -- (-1,-1);
\draw (C) -- (0,0);
\draw (-1,-1) -- (0,0);

\end{scope}


\begin{scope}[shift={(0,-1)}]

\draw[->] (-1,-1) -- (1,-1) node[midway,above]{Rotation at $X$};
\draw[->] (1,-2) -- (-1,-2) node[midway,above]{Rotation at $Y$};

\end{scope}


\begin{scope}[yscale=2,shift={(3.5,0)}]

\node (Y) at (1,-1) {$\bullet$};
\node at (Y) [right = 1mm of Y]{$Y$};
\node (X) at (0,0) {$\bullet$};
\node at (X) [right = 1mm of X]{$X$};
\node [regular polygon,regular polygon sides=3,draw] (A) at (-1,-1) {$A$};
\node [regular polygon,regular polygon sides=3,draw] (B) at (0,-2) {$B$};
\node [regular polygon,regular polygon sides=3,draw] (C) at (2,-2) {$C$};

\draw (A) -- (0,0);
\draw (B) -- (1,-1);
\draw (C) -- (1,-1);
\draw (1,-1) -- (0,0);

\end{scope}

\end{tikzpicture}
\]
\end{figure}

\emph{Rotation} of binary trees is the operation shown in Figure~\ref{fig:rotation}. Here $X$ and $Y$ are nodes and $A$, $B$, and $C$ represent subtrees. In the left-hand tree, $A$ is the left subtree of $X$, $B$ is the right subtree of $X$, and $C$ is the right subtree of $Y$; in the right-hand tree, $A$ is still the left subtree of $X$ and $C$ is still the right subtree of $Y$, but $B$ is now the left subtree of $Y$. Note that the trees depicted in this figure might be subtrees of larger trees; that is, rotation can be done as an operation on subtrees of larger trees.

Given the situation of Figure~\ref{fig:rotation}, noting that the trees may be subtrees of larger trees, we call the rotation at $X$ a \emph{forwards rotation} and the rotation at $Y$ a \emph{backwards rotation}. The \emph{Tamari lattice} $T_{n}$ is the poset of binary trees of size $n$ with covering relations given by forwards rotation.

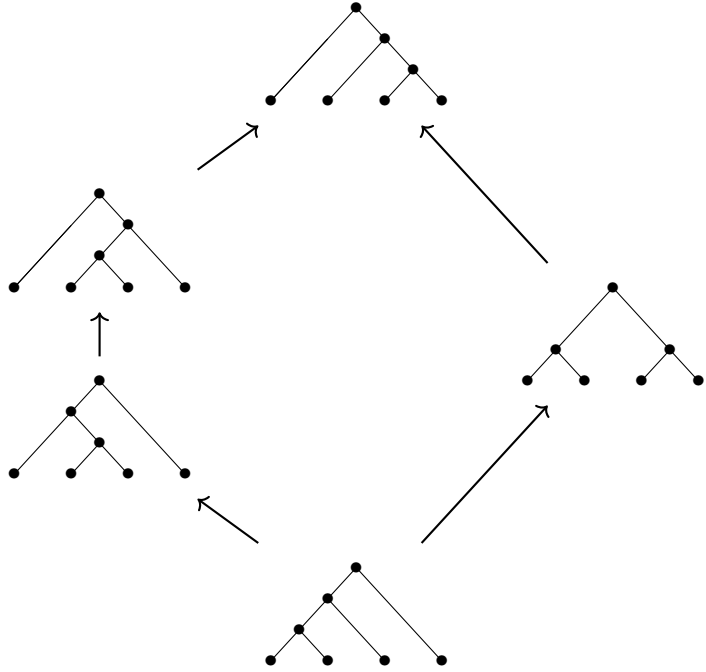
\begin{figure}
\caption{The Tamari lattice $T_{4}$}\label{fig:tamari}
\[
\begin{tikzpicture}[xscale=2,scale=0.75,yscale=1.1]


\begin{scope}[xscale=0.5,scale=0.5]

\draw (0,0) -- (-1,-1);
\draw (0,0) -- (3,-3);
\draw (-1,-1) -- (-2,-2);
\draw (-1,-1) -- (1,-3);
\draw (-2,-2) -- (-1,-3);
\draw (-2,-2) -- (-3,-3);

\node (x) at (0,0) {$\bullet$};
\node at (-1,-1) {$\bullet$};
\node (y) at (3,-3) {$\bullet$};
\node at (-2,-2) {$\bullet$};
\node at (1,-3) {$\bullet$};
\node at (-1,-3) {$\bullet$};
\node (z) at (-3,-3) {$\bullet$};

\node[draw=none,fit=(x) (y) (z)](b){};

\end{scope}


\begin{scope}[xscale=0.5,scale=0.5,shift={(9,9)}]

\draw (0,0) -- (-2,-2);
\draw (0,0) -- (3,-3);
\draw (2,-2) -- (1,-3);
\draw (-2,-2) -- (-1,-3);
\draw (-2,-2) -- (-3,-3);

\node (x) at (0,0) {$\bullet$};
\node at (2,-2) {$\bullet$};
\node (y) at (3,-3) {$\bullet$};
\node at (-2,-2) {$\bullet$};
\node at (1,-3) {$\bullet$};
\node at (-1,-3) {$\bullet$};
\node (z) at (-3,-3) {$\bullet$};

\node[draw=none,fit=(x) (y) (z)](l1){};

\end{scope}


\begin{scope}[xscale=0.5,scale=0.5,shift={(-9,6)}]

\draw (0,0) -- (-1,-1);
\draw (0,0) -- (3,-3);
\draw (-1,-1) -- (-3,-3);
\draw (-1,-1) -- (1,-3);
\draw (0,-2) -- (-1,-3);

\node (x) at (0,0) {$\bullet$};
\node at (-1,-1) {$\bullet$};
\node (y) at (3,-3) {$\bullet$};
\node at (0,-2) {$\bullet$};
\node at (1,-3) {$\bullet$};
\node at (-1,-3) {$\bullet$};
\node (z) at (-3,-3) {$\bullet$};

\node[draw=none,fit=(x) (y) (z)](r1){};

\end{scope}


\begin{scope}[xscale=0.5,scale=0.5,shift={(-9,12)}]

\draw (0,0) -- (-3,-3);
\draw (0,0) -- (3,-3);
\draw (-1,-1) -- (-3,-3);
\draw (1,-1) -- (0,-2);
\draw (0,-2) -- (1,-3);
\draw (0,-2) -- (-1,-3);

\node (x) at (0,0) {$\bullet$};
\node at (1,-1) {$\bullet$};
\node (y) at (3,-3) {$\bullet$};
\node at (0,-2) {$\bullet$};
\node at (1,-3) {$\bullet$};
\node at (-1,-3) {$\bullet$};
\node (z) at (-3,-3) {$\bullet$};

\node[draw=none,fit=(x) (y) (z)](r2){};

\end{scope}


\begin{scope}[xscale=0.5,scale=0.5,shift={(0,18)}]

\draw (0,0) -- (-3,-3);
\draw (0,0) -- (3,-3);
\draw (-1,-1) -- (-3,-3);
\draw (1,-1) -- (-1,-3);
\draw (2,-2) -- (1,-3);

\node (x) at (0,0) {$\bullet$};
\node at (1,-1) {$\bullet$};
\node (y) at (3,-3) {$\bullet$};
\node at (2,-2) {$\bullet$};
\node at (1,-3) {$\bullet$};
\node at (-1,-3) {$\bullet$};
\node (z) at (-3,-3) {$\bullet$};

\node[draw=none,fit=(x) (y) (z)](t){};

\end{scope}


\draw[->,thick] (b) -- (l1);
\draw[->,thick] (b) -- (r1);
\draw[->,thick] (l1) -- (t);
\draw[->,thick] (r1) -- (r2);
\draw[->,thick] (r2) -- (t);

\end{tikzpicture}
\]
\end{figure}

We now define the weak Bruhat order on the symmetric group $S_{n}$. A more general definition of the weak Bruhat order which applies to all Coxeter groups will be defined in Section~\ref{sect:ex:camb}. Given a permutation $\sigma \in S_{n}$, an \emph{inversion} of $\sigma$ is a pair $\{i, j\} \in \binom{[n]}{2}$ with $i < j$, such that $\sigma(j) < \sigma(i)$. The set of all inversions of $\sigma$ is denoted $I(\sigma)$. The \emph{weak Bruhat order} on the symmetric group is the partial order defined by $\sigma \leqslant \tau$ if and only if $I(\sigma) \subseteq I(\tau)$.

The map $\psi\colon S_{n} \to T_{n}$ from the weak Bruhat order on permutations to the Tamari lattice on binary trees is defined recursively. We follow the exposition in \cite{thomas-bst}. For a string of $p$ distinct numbers $a_{1}\dots a_{p}$, we write $\std{a_{1}\dots a_{p}}$ for the string consisting of $[p]$ written in the same order as $a_{1}, \dots, a_{p}$. We use one-line notation for permutations $\sigma \in S_{n}$, so that $\sigma = w_{1}w_{2}\dots w_{n}$ means that $\sigma(i) = w_{i}$. For $n = 0$, $\psi$ applied to the empty permutation gives the empty tree. Then, for $n \geqslant 1$, given a permutation $\sigma = a_{1} \dots a_{p} n b_{1} \dots b_{q}$, we define $\psi(\sigma)$ to be the binary tree where the root node has left subtree $\psi(\std{a_{1} \dots a_{p}})$ and right subtree $\psi(\std{b_{1} \dots b_{q}})$. This map is order-preserving, since adding an inversion to the permutation either corresponds to a forwards rotation of the tree, or has no effect. For example, going from $12$ to $21$ gives a forwards rotation of the tree, whereas going from $132$ to $231$ has no effect.

\begin{theorem}[\cite{reading_cambrian,bw_shell_2}]
The map $\psi\colon S_{n} \to T_{n}$ is a quotient map by a lattice congruence.
\end{theorem}

\subsection{Cambrian lattices}\label{sect:ex:camb}

We now examine a generalisation of the map from Section~\ref{sect:ex:perms->trees} due to Reading \cite{reading_cambrian}. A \emph{Coxeter group} is a group $W$ defined by a set of generators $S$ and relations $(st)^{m(s, t)} = 1$ for $t, s \in S$, where $m(s, t) = 1$ for $s = t$ and $m(s, t) \in [2, \infty]$ otherwise. The elements of $S$ are called \emph{simple reflections} and conjugates of simple reflections are called \emph{reflections}. Important examples of Coxeter groups include Weyl groups and finite reflection groups. More details can be found in \cite{bb_coxeter}.

An element of $W$ can be written in several different ways as a word with letters in $S$. A word $a$ for an element $w \in W$ is called \emph{reduced} if it has a minimal number of letters amongst all words representing $w$. The length of a reduced word for $w$ is called the \emph{length} $\ell(w)$ of $w$. A finite Coxeter group has an unique element of maximal length, which is referred to as the \emph{longest element} and denoted $w_{0}$.

We consider a partial order on Coxeter groups, known as the \emph{weak Bruhat order}, generalising the order on the symmetric group considered in Section~\ref{sect:ex:perms->trees}. There is a family of lattice quotients of this partial order known as \emph{Cambrian lattices} \cite{reading_cambrian}. The prototype for this family of quotients is the map from permutations to binary trees given in Section~\ref{sect:ex:perms->trees}.

We again let $W$ be a Coxeter group with $S$ the set of simple reflections. Letting $w \in W$, the \emph{left inversion set} $I(w)$ of $w$ is defined to be the set of all reflections $t$ such that $\ell(tw) < \ell(w)$. We have that $|I(w)| = \ell(w)$. The \emph{left descent set} of $w$ is $I(w) \cap S$. There are an analogous right inversion set and right descent set: the \emph{right descents} of $w$ are the simple reflections $s \in S$ such that $\ell(ws) < \ell(w)$.

The \emph{right weak Bruhat order} on $W$ is the partial order on $W$ where $v \leqslant w$ if and only if $I(v) \subseteq I(w)$. Note that, unfortunately, the \textit{right} weak Bruhat order is defined in terms of \textit{left} inversion sets. Equivalently, the right weak Bruhat order is the partial order with covering relations $ws \lessdot w$ for every right descent $s$ of $w$. Again equivalently, $v \leqslant w$ in the right weak Bruhat order if and only if there is a reduced word for $v$ which is a prefix of a reduced word for $w$. There is an analogous left weak Bruhat order, which is isomorphic to the right weak Bruhat order via the map $w \mapsto w^{-1}$. Henceforth, when we say `weak Bruhat order', we will mean the right weak Bruhat order.

The \emph{Coxeter diagram} of $W$ with respect to $S$ is the graph whose vertex set is $S$ and whose edges $\{s, t\}$ are given by the pairs such that $m(s, t) \geqslant 3$. The edges $\{s, t\}$ with $m(s, t) \geqslant 4$ are labelled by the number $m(s, t)$. An \emph{orientation} of a Coxeter diagram $G$ is a directed graph $G^{\to}$ with the same vertex set as $G$ with one directed edge for each edge of $G$. The Coxeter diagram corresponding to the symmetric group $S_{n + 1}$ is the $A_{n}$ Dynkin diagram.

We are now in a position to define the Cambrian congruences from \cite{reading_cambrian}. Let $W$ now be a finite Coxeter group with Coxeter diagram~$G$. Further, let $G^{\to}$ be an orientation of $G$. The \emph{Cambrian congruence} $\Theta(G^{\to})$ is the smallest lattice congruence on the weak Bruhat order of $W$ such that, for a directed edge $s \to t$ in $G^{\to}$ with label $m(s, t)$, $t$ is equivalent to \[\underbrace{tsts\dots}_{m(s, t) - 1}.\] The associated \emph{Cambrian lattice} is the quotient $W/\Theta(G^{\to})$. In the case of the example from Section~\ref{sect:ex:perms->trees}, the orientation of the $A_{n}$ Dynkin diagram giving the Tamari lattice is the one where all arrows point in the same direction.

\begin{figure}
\caption{A Cambrian congruence on the weak Bruhat order on $S_{2}$}\label{fig:ex:camb}
\[
\begin{tikzpicture}[
	frames/.style args = {#1/#2}{minimum height=#1,
               minimum width=#2+\pgfkeysvalueof{/pgf/minimum height},
               draw, rounded corners=3mm, fill=yellow, opacity=0.6,
               sloped},
               ]
               
\begin{scope}[shift={(-2,0)}]

	\begin{pgfonlayer}{foreground}
\node(0) at (0,0) {$e$};
\node(a) at (-1,1) {$s_{1}$};
\node(b) at (1,1) {$s_{2}$};
\node(c) at (-1,2) {$s_{1}s_{2}$};
\node(d) at (1,2) {$s_{2}s_{1}$};
\node(1) at (0,3) {$w_{0}$};

\draw[->] (0) -- (a);
\draw[->] (0) -- (b);
\draw[->] (a) -- (c);
\draw[->] (b) -- (d);
\draw[->] (c) -- (1);
\draw[->] (d) -- (1);
	\end{pgfonlayer}

    \begin{pgfonlayer}{main}
\path   let \p1 = ($(b.center)-(d.center)$),
            \n1 = {veclen(\y1,\x1)} in
            (b) -- 
            node[frames=7mm/\n1] {} (d);
    \end{pgfonlayer}

\end{scope}


\begin{scope}[shift={(2,0)}]

\node(0) at (0,0) {$e$};
\node(a) at (-1,1) {$s_{1}$};
\node(bd) at (1,1.5) {$\{s_{2}, s_{2}s_{1}\}$};
\node(c) at (-1,2) {$s_{1}s_{2}$};
\node(1) at (0,3) {$w_{0}$};

\draw[->] (0) -- (a);
\draw[->] (0) -- (bd);
\draw[->] (a) -- (c);
\draw[->] (c) -- (1);
\draw[->] (bd) -- (1);

\end{scope}

\end{tikzpicture}
\]
\end{figure}
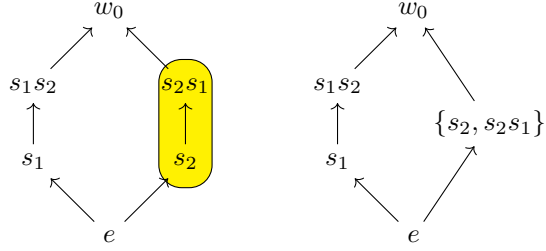

\begin{example}\label{ex:camb}
We give an example of a Cambrian congruence and the resulting lattice. Consider the orientation \[1 \to 2\] of the $A_{2}$ Dynkin diagram. If we let $s_{1}$ and $s_{2}$ be the simple reflections corresponding to the relevant vertices. We have that $m(s_{1}, s_{2}) = 3$, so under the Cambrian congruence we have that $s_{2}$ is equivalent to $s_{2}s_{1}$. This congruence and resulting lattice is shown in Figure~\ref{fig:ex:camb}. Note that the Cambrian lattice is the Tamari lattice from Figure~\ref{fig:tamari}.
\end{example}

Cambrian lattices are sublattices of the weak order, which is not generally true for lattice quotients, as we know from Example~\ref{ex:subpos_not_sublat}.

\begin{theorem}[{\cite[Theorem~1.1, Theorem~1.2]{reading_sortable}}]
Let $G$ be a Coxeter diagram with $W$ the Coxeter group of $G$. Let $\Theta(G^{\to})$ be a Cambrian congruence on $W$ for some orientation $G^{\to}$ of $G$. Then the Cambrian congruence $\Theta(G^{\to})$ is a lattice congruence on $W$ and the Cambrian lattice $W/\Theta(G^{\to})$ is a sublattice of $W$.
\end{theorem}

See also the related paper \cite{reading_weak}. Cambrian quotients of infinite Coxeter groups were considered in \cite{rs_seicg}.

\begin{remark}
Cambrian lattices and quotients can be interpreted in the representation theory of algebras in \cite{it,dirrt,irrt,g-mc,gyoda_nsl}. See also \cite{mizuno-preproj}, which realises the weak Bruhat order in terms of the representation theory of preprojective algebras.
\end{remark}

\subsection{Higher-dimensional Cambrian maps}

The map $\psi$ from Section~\ref{sect:ex:perms->trees} also extends to higher dimensions. The higher Bruhat orders $\bruhat{n}{d}$ are a family of partial orders such that $\bruhat{n}{1}$ is the weak Bruhat order on $S_{n}$. These were introduced in \cite{ms}. Similarly, the higher Stasheff--Tamari orders $\stash{n}{d}$ are a family of partial orders which coincide with the Tamari lattice $T_{n}$ in the case $d = 2$. These orders were defined in \cite{kv-poly,er} and a good introduction can be found in \cite{rr}. There were originally two different descriptions of these orders, but these were shown to be the same in \cite{njw-equal}. In \cite{kv-poly}, Kapranov and Voevodsky defined an order-preserving map $f \colon \bruhat{n}{d} \to \stash{n + 2}{d + 1}$, which is a higher-dimensional version of the map $\psi$ from Section~\ref{sect:ex:perms->trees}; that is, $f = \psi$ for $d = 1$. This map was studied in \cite{rambau,thomas-bst}. The following conjecture is open.

\begin{conjecture}[{\cite[Theorem~4.10]{kv-poly}}]
The map \[f\colon \bruhat{n}{d} \to \stash{n + 2}{d + 1}\] is a quotient by a weak order congruence.
\end{conjecture}

\begin{remark}
The higher Bruhat orders $\bruhat{n}{d}$ and higher Stasheff--Tamari orders $\stash{n + 2}{d + 1}$ are not only posets but $n - d$ categories, with the order relations giving the one-dimensional morphisms. Being a quotient by a weak order congruence effectively means that the map $f$ is surjective on the order relations. But the map $f$ should in fact be surjective on morphisms of all dimensions in $\stash{n + 2}{d + 1}$. However, showing surjectivity on the relations ought to suffice, since the higher-dimensional morphisms in $\bruhat{n}{d}$ and $\stash{n + 2}{d + 1}$ correspond to relations in these posets for larger values of $d$.
\end{remark}

It is known that $f$ cannot be a quotient by an order congruence due to the following fact.

\begin{proposition}[{\cite[Section~6]{thomas-bst}}]
The fibres of the map $f$ are not always intervals.
\end{proposition}

There is also a map $g \colon \bruhat{n}{d + 1} \to \stash{n}{d}$ which factors through the dual of the map $f$ \cite[Proposition~7.1]{thomas-bst}, \cite[Remark~34]{njw-hbo}, for which the following result holds.

\begin{theorem}[{\cite[Theorem~5.3]{njw-hbo}}]
The map $g \colon \bruhat{n}{d + 1} \to \stash{n}{d}$ is a quotient by a weak order congruence.
\end{theorem}

The surjectivity of the map $g$ is originally due to \cite[Theorem~3.5]{rs-baues}. This quotient map also appears in the context of integrable systems in \cite{dm-h}, who use it to give a definition of the higher Stasheff--Tamari orders as a quotient of the higher Bruhat orders.

\subsection{Type \texorpdfstring{$B$}{B} weak order}

In \cite{simion}, Simion gives a quotient of the weak order in type $B$ which is isomorphic to the weak order in type~$A$, which we now explain. Let $W_{B_{n}}$ be the Coxeter group with Coxeter diagram $B_{n}$. This Coxeter group is the group of \emph{signed permutations}, meaning permutations $x_{-n} \dots x_{-1}x_{1} \dots x_{n}$ of $\pm [n] := \{\pm 1, \dots, \pm n\}$ such that $x_{-i} = - x_{i}$ for every $i \in [n]$. Signed permutations can therefore be represented in \emph{abbreviated one-line notation} by $x_{1} \dots x_{n}$. We write $|x_{i}|$ for the absolute value of $x_{i}$ in the natural way: $|j| = j$ and $|-j| = j$ for $j \in [n]$.

Given $\sigma = \sigma_{1} \dots \sigma_{n} \in W_{B_{n}}$, let $\sigma^{+}$ consist of the subword of $\sigma$ consisting of positive entries and $\sigma^{-}$ consist of the subword of $\sigma$ consisting of negative entries. Furthermore, we let $\alpha_{\sigma} = \sigma^{+}\sigma^{-}$ and $\beta_{\sigma} = \sigma^{-}\sigma^{+}$.

\begin{lemma}[{\cite[Lemma~3]{simion}}]\label{lem:simion_partition}
The set of intervals \[\spt{[\alpha_{\sigma}, \beta_{\sigma}] \st \sigma \in W_{B_{n}}}\] forms a partition of $W_{B_{n}}$.
\end{lemma}

If $\equ$ is the equivalence relation given by this partition, then we have the following result.

\begin{proposition}[{\cite[Proposition~4]{simion}}]
The equivalence relation $\equ$ is an order congruence on $W_{B_{n}}$ and the quotient poset $W_{B_{n}}/{\equ}$ is isomorphic to $W_{A_{n}}$, the weak Bruhat order on the symmetric group.
\end{proposition}

There is a nice description of the intervals in the partition.

\begin{proposition}[{\cite[Observation~2]{simion}}]
Let $\sigma \in W_{B_{n}}$ be such that $\sigma^{+}$ has $k$ letters. Then $[\alpha, \beta] \cong L(n - k, k)$.
\end{proposition}

Here $L(n - k, k)$ is the poset of Young diagrams from Section~\ref{sect:ex:young}.

\subsection{Strong Bruhat order}

We now consider a family of quotients of the strong Bruhat order on a Coxeter group $W$. Given a Coxeter group $W$, the \emph{strong Bruhat order} on $W$ is defined as follows. We have that $u \leqslant w$ if some reduced word of $w$ contains a reduced word for $u$ as a subword. In fact, in this case, every reduced word for $w$ will contain a reduced word for~$u$. Unlike the weak Bruhat order, the strong Bruhat order is not always a lattice, and even not always either a meet-semilattice or a join-semilattice \cite[Section~3.2]{bb_coxeter}.

For any subset $J \subset S$, the subgroup of $W$ generated by $J$ is another Coxeter group, known as a \emph{parabolic subgroup} and denoted $W_{J}$. The following proposition shows that the quotient of the strong Bruhat order on a finite Coxeter group using the two-sided cosets from two parabolic subgroups is an order quotient.

\begin{proposition}[{\cite[Proposition~31]{reading_order}}]
For any $w \in W$ and $J, K \subseteq S$, the double cosets $W_{J}wW_{K}$ form an order congruence of the strong Bruhat order on~$W$.
\end{proposition}

For finite Coxeter groups, quotients of the strong Bruhat order by parabolic subgroups were classified in \cite{newton}.

\subsection{Bruhat interval polytopes}

We finally consider another family of quotients of the weak Bruhat orders, which nevertheless use the strong Bruhat order in their definition. For a permutation $w \in S_{n}$, the associated \emph{Bruhat interval polytope}  \cite[Definition~6.9]{kw_kt} is defined as the convex hull \[\mathfrak{Q}_{w} := \conv\spt{\mathbf{u} := (u^{-1}(1), \dots, u^{-1}(n)) \st u \leqslant w} \subset \mathbb{R}^n,\] where here $u \leqslant w$ means with respect to the strong Bruhat order on $S_{n}$.
Recall that the \emph{convex hull} of a set of points is the smallest convex set containing them, where a set is \emph{convex} if it contains the straight line segment between any two of its points;
a \emph{convex polytope} is the convex hull of a finite set of points.
By looking at the edges of the Bruhat interval polytope, one obtains the following partial order on the set $[e, w]$.

\begin{definition}{\cite[Definition~1.1]{gaetz_bip}}
The poset $(P_{w}, \leqslant_{w})$ has the interval $[e, w]$ in the strong Bruhat order as its underlying set, with covering relations $u \lessdot_{w} v$ whenever $\mathfrak{Q}_{w}$ has an edge between the vertices $\mathbf{u}$ and $\mathbf{v}$ and $\ell(u) < \ell(v)$.
\end{definition}

These posets are lattices.

\begin{theorem}[{\cite[Theorem~4.5]{gaetz_bip}}]
For all $w \in S_{n}$, we have that $P_{w}$ is a lattice.
\end{theorem}

One can realise $P_{w}$ as a certain quotient poset of the weak Bruhat order on $S_{n}$.
In order to do this, we will need the following notion.
The \emph{normal fan} $N(\mathfrak{Q})$ of a polytope $\mathfrak{Q} \subset \mathbb{R}^n$, is the set of cones $\spt{C(F) \st F \text{ a face of } \mathfrak{Q}}$, where \[C(F) := \spt{\mathbf{x} \in \mathbb{R}^n \st F \subseteq \argmax_{\mathbf{x}' \in \mathfrak{Q}} \ip{\mathbf{x}}{\mathbf{x}'}},\] and $\argmax_{\mathbf{x}' \in \mathfrak{Q}} \ip{\mathbf{x}}{\mathbf{x}'} := \spt{\mathbf{x}' \in \mathfrak{Q} \st \ip{\mathbf{x}}{\mathbf{x}'} \text{ is maximal}}$ \cite[Ex.~7.3]{ziegler}.

Elements $u$ of $P_{w}$ correspond to vertices $\mathbf{u}$ of $\mathfrak{Q}_{w}$, which correspond to maximal cones $C(\mathbf{u})$ in $N(\mathfrak{Q}_{w})$.
We have that $\mathfrak{Q}_{w_{0}}$ is the permutohedron, whose normal fan $N(\mathfrak{Q}_{w_{0}})$ is the \emph{braid hyperplane arrangement}, which has defining hyperplanes $x_{i} - x_{j} = 0$.
In the normal fan of $\mathfrak{Q}_{w}$, each cone of $N(\mathfrak{Q}_{w})$ is a union of cones of $N(\mathfrak{Q}_{w_{0}})$.
Hence, we obtain an equivalence relation $\theta_{w}$ on the weak Bruhat order $P_{w_{0}}$, in which the equivalence class of $u \in P_{w_{0}}$ is the set of $v \in P_{w_{0}}$ whose corresponding cones $C(\mathbf{v})$ lie in the same cone of~$N(\mathfrak{Q}_{w})$ as $C(\mathbf{u})$.
The elements of $P_{w}$ are then in bijection with the $\theta_{w}$-equivalence classes of~$P_{w_{0}}$.

\begin{theorem}[{\cite[Theorem~4.3, Theorem~4.5]{gaetz_bip}}]
We have that $P_{w} \cong P_{w_{0}}/\theta_{w}$.

Moreover, for $u, v \in P_{w_{0}}$, we have \[[u]_{\theta_{w}} \join [v]_{\theta_{w}} = [u \join v]_{\theta_{w}}.\]
\end{theorem}

Hence $\theta_{w}$ is a join-semilattice congruence, and so is a III-congruence on $\overline{P_{w_{0}}}$.

\section{Comparison}\label{sect:comparison}

We finally compare the differing notions of poset congruence that we have surveyed. Different types of poset congruence have different properties, and some are stronger than others.

In Table~\ref{table}, we compare the different properties of the notions of congruence.
\begin{itemize}
\item \textbf{Self-dual?} A tick indicates that $\theta$ is a congruence on $P$ if and only if it is a congruence on $\dual{P}$.
\item \textbf{Preserves grading?} We have that $\quot$ is a graded poset if $P$ is a connected graded poset.
\item \textbf{Quotient map strong?} The canonical map $P \to \quot$ is strong. Equivalently, one does not need to take the transitive closure $\trans{\qrel}$ of the quotient relation $\qrel$ on $\quot$ to define the quotient relation.
\item \textbf{Closed under $\cap$?} The intersection $\equ_{1} \cap \equ_{2}$ is a congruence if $\equ_{1}$ and $\equ_{2}$ are congruences.
\item \textbf{Infinite posets?} If $\equ$ is a congruence on an infinite poset $P$, then $\quot$ is still a well-defined poset.
\item \textbf{Lattice congruences?} If $\equ$ is a congruence on a lattice $L$, then $\equ$ is a lattice congruence on $L$.
\item \textbf{Other requirements.} Other conditions needed on the poset $P$.
\end{itemize}
References for some of these facts are as follows, corresponding to the superscripts.
\begin{enumerate}
\item \cite[Theorem~6.3]{bj_res}.\label{ref:comp:int}
\item \cite[Lemma~11]{halas}.\label{ref:wsc:int}
\item \cite[p.197]{halas}.\label{ref:ord:int}
\item \cite[Theorem~4.8]{hl_vop}.\label{ref:hl:int}
\item \cite[p.347]{hl_vop}.\label{ref:hl:lat}
\item \cite[2.6]{kolibiar}.\label{ref:kol:full}
\item \cite[2.1]{kolibiar}.\label{ref:kol:lat}
\end{enumerate}
For reasons of space, the remaining facts are left as exercises.

In Figure~\ref{fig:implications}, we display the implications that hold between different notions of poset congruence. Here `Haviar--Lihov\'a congruence $\implies$ Order congruence' means that every Haviar--Lihov\'a congruence is an order congruence, and so on. We label the arrows with the additional assumptions that are needed, if any.
\begin{itemize}
\item Strong: the implication holds if the quotient map $P \to \quot$ is strong.
\item $\{\hat{0}\}$: the implication holds if the poset has a unique minimal element $\hat{0}$ which is in its own equivalence class $\{\hat{0}\}$. Note that in the case of the orbits of a group of automorphism, if there exists a unique minimal element $\hat{0}$, then it is automatic that it is in its own equivalence class $\{\hat{0}\}$.
\item Finite: the implication holds if the poset is finite.
\end{itemize}
We provide a citation if the result can be found explicitly found in the literature. The remaining implications are left as exercises.

A good illustration of the fact that several notions of congruence are needed is illustrated by the following proposition.

\begin{proposition}\label{prop:orb_int_ord=triv}
Let $P$ be a poset with $\theta$ an order congruence which is also the set of orbits of an automorphism group. Then $\theta$ is the identity relation.
\end{proposition}

Order congruences and orbits of groups of automorphisms are both natural classes of equivalence relations by which to quotient a poset, but the two only intersect in the identity relation. Hence both notions are needed separately if the spectrum of different types of congruence that arise in examples are to be captured.

\printbibliography

\end{document}